\documentclass{article}

\usepackage{amsmath,amssymb,amsbsy,}
\usepackage{amsthm}
\usepackage[mathscr]{eucal}

\newtheorem{theorem}{Theorem}\numberwithin{theorem}{section}
\newtheorem{lem}[theorem]{Lemma}

\newtheorem{exc}[theorem]{Exercise}
\newtheorem{rem}[theorem]{Remark}

\newcommand{\Aff}{\mathbb A}
\newcommand{\C}{\mathbb C}
\newcommand{\G}{\mathbb G}
\newcommand{\T}{\mathbb T}
\newcommand{\Q}{\mathbb Q}
\newcommand{\Z}{\mathbb Z}
\newcommand{\M}{\mathbb M}
\newcommand{\PP}{\mathbb P}
\newcommand{\simby}[1]{\buildrel{#1}\over\sim}    % for eqce relns
\newcommand{\simT}{\simby{\T}}
\newcommand{\tensor}{\otimes}
\newcommand{\onto}{\twoheadrightarrow}
\newcommand{\Span}[1]{\left<#1\right>}

\newcommand{\dd}{\mathrm d}
\newcommand{\om}{\omega}

\newcommand{\al}{\alpha}
\newcommand{\be}{\beta}
\newcommand{\la}{\lambda}
\newcommand{\si}{\sigma}
\newcommand{\bmu}{\boldsymbol\mu}
\newcommand{\checku}{\check{u}}
\newcommand{\checkv}{\check{v}}
\newcommand{\Gm}{\mathbb G_m}
\DeclareMathOperator{\GL}{GL}
\DeclareMathOperator{\Gr}{Gr}
\DeclareMathOperator{\aGr}{aGr}
\DeclareMathOperator{\Pf}{Pf}
\DeclareMathOperator{\SL}{SL}
\DeclareMathOperator{\SO}{SO}
\DeclareMathOperator{\Sym}{Sym}

\newenvironment{mycase}[1]{\medskip\paragraph{{\bf Case #1.}\enspace}}{\medskip}

\newenvironment{pf}{\paragraph{{\bf Proof.}\enspace}}{{\hfill Q.E.D.}\medskip}
\numberwithin{figure}{section}
\numberwithin{equation}{section}

\title{Diptych varieties. II: Apolar varieties}
\author{Gavin Brown and Miles Reid}
\date{July 2015}

\begin{document}

\maketitle

\begin{abstract}
This paper constructs all the diptych varieties with $de\le4$ (see
\cite{BR1}, Main Theorem~3.3). Our construction involves several new
classes of Gorenstein almost homogeneous spaces for
$\GL(2)\times\G_m^r$, in particular two infinite series arising from the
algebra of apolarity.
\end{abstract}

%\keywords{Keywords: Almost homogeneous space, Gorenstein variety, unprojection, apolar geometry, Mori flip.}
%\\ 2010 {\em Mathematics Subject Classification}: 14M17 (Primary);  14E99 (Secondary).}

\setcounter{tocdepth}{2}
\tableofcontents

\subsection*{Diptych varieties and Mori flips}

We introduced {\em diptych varieties} in \cite{BR1}, motivated by our
attempts to understand Mori's explicit calculations \cite{M} in the
Picard group of a 3-fold extremal neighbourhood. Mori's argument
associates a 2-step continued fraction expansion $[d,e,d,\dots]$ with an
extremal neighbourhood. Roughly, for $C=\PP^1\subset X$ a flipping curve
of Type~A in a 3-fold $X$ with two terminal singularities $P,Q\in C$ of
type $cA_n/\bmu_r$ and a pair of divisors transverse to $C$ at $P$ and
$Q$ respectively, Mori sets up a `continued division' algorithm that
constructs a sequence of divisors $F_{2i}\sim F_{2i-1}-dF_{2i-2}$,
$F_{2i+1}\sim F_{2i}-eF_{2i-1}$, and proves that it terminates in the
set theoretic equality $C=F_k\cap F_{k+1}$ for some $k$. This expresses
a flipping curve $C$ as the base locus of a pencil of divisors, and
hence proves the existence of the flip of $C\subset X$, showing moreover
that it can in principle be computed as the normalisation of the pencil.
Diptych varieties are {\em key varieties} for the $\G_m$ cover of these
{\em Type~A flips}: flips arise as regular pullbacks from diptychs after
some massaging; see \cite{Ki} \S11 (especially 11.2) and \cite{BR4} for
details of this last step from diptychs to extremal neighbourhoods.

For completeness, we give some details in \S\ref{s!dip1} of what we
understand by a diptych variety; in brief, each is an affine 6-fold
$V_{ABLM}$ arising as a 4-parameter deformation of a {\em tent}, a
reducible Gorenstein toric surface consisting of a cycle $T=S_0\cup
S_1\cup S_2\cup S_3$ of four affine toric components meeting along their
1-dimensional strata; the four deformation parameters smooth the axes of
transverse intersections of the cycle. A diptych variety is
characterised by three natural numbers $d,e,k$, or by a 2-step recurrent
continued fraction $[d,e,d,\dots]$ to $k$ terms -- of course, these
correspond to the $d,e,k$ of Mori's continued division algorithm.

Theorem~1.1 of \cite{BR1} asserts that a diptych variety exists for any
$d,e,k$ (with the bounds of \cite{BR1}, Theorem~3.3, (3.7) on $k$ in the
cases $de\le3$). In the main case $de>4$ and $d,e\ge2$, we proved this
in \cite{BR1}, Section~5. In \cite{BR3} we treat the cases $de>4$ with
$d$ or $e=1$ using variants of the same methods. This paper constructs
diptych varieties in the remaining cases $de\le4$, fulfilling the
promise of \cite{BR1}, Theorem~1.1, and providing key varieties for the
remaining extremal neighbourhoods of Type~A.

\subsection*{Apolar geometry}

The diptych varieties with $de=4$ have a beautiful description in terms
of key 5-folds $V_k\subset\Aff^{k+5}$ that play a principal role in
this paper (see \S\ref{s!1}, and especially \ref{ss!Vk}). These are {\em
almost homogeneous spaces} that are easy to describe based on the
algebra of apolarity, and we offer several alternative approaches. With
a final unprojection argument, any of these descriptions is enough to
prove the existence of diptych varieties with $de=4$.

Geometrically, the $V_k$ are almost homogeneous spaces for the group
$G=\GL(2)\times\G_m$: each is the closure of the orbit of an `apolar'
vector in a reducible representation of $G$, and we refer to them as
{\em apolar varieties}, as yet with no general formal definition, but
see \ref{ss!Vk}. It would be interesting to know whether apolar
varieties such as the $V_k$ and the $W_d$ introduced in \ref{ss!32}
arise naturally in other parts of geometry and representation theory; we
see similar phenomena in other calculations in codimension $\ge4$, and
this type of {\em apolar geometry} should apply more widely.

From the point of view of equations, we express the $V_k$ using a
generalised form of Cramer's rule. This provides all the equations of
$V_k$ in closed form, in contrast to the small subset of Pfaffian
equations that we get away with in \cite{BR1}. The varieties $V_k$ are
serial unprojections, although this does not itself provide all the
equations directly.

\S\ref{s!lt3} introduces a second series of apolar varieties, this time
almost homo\-geneous \hbox{7-folds} $W_d\subset\Aff^{d+9}$, and applies
them as models for dip\-tychs with $k=2$. With a single additional
unprojection, they also provide a format for diptychs with $k=3$
involving {\em crazy Pfaffians}, reminiscent of Riemenschneider's
`quasi-determinants' \cite{R}; see \ref{ss!33} where we discuss the
equations in terms of {\em floating factors}. \S\ref{s!k45} handles the
few remaining cases with $k=4,5$ and $de=3$, where unprojection methods
and pentagrams provide the equations directly. Rather than our apolar
varieties $V_k$ and $W_d$ given by serial unprojection, these cases are
most naturally described as regular pullbacks from a parallel
unprojection key variety, a 10-fold $W\subset\Aff^{16}$.

\subsection*{Gorenstein rings in high codimension}

Gorenstein rings arise naturally in geometry as homogeneous coordinate
rings of Fanos, Calabi--Yaus, regular canonical $n$-folds, and other
constructions -- and, most notably for our purposes here, of 3-fold
extremal neighbourhoods. Thus a supply of model Gorenstein rings, with
explicit information about their generators and relations, gradings and
so on, is of practical importance. It is hard to construct Gorenstein
rings in high codimension in general; there is no practical
classification beyond codimension~3 (although see \cite{Reid2,Reid3} for
a first structure theorem in codimension~4). Grojnowski and Corti and
Reid \cite{CR} study weighted homo\-geneous spaces or closed orbits in
highest weight representations of semisimple algebraic groups, in
particular for $\SL(5)$ and $\SO(10)$; Qureshi and Szendr{\H o}i
\cite{QS} generalise these to more classes of examples. The almost
homogeneous spaces $V_k$ in \S\ref{s!1} (dimension~5, codimension~$k$),
$W_d$ in \S\ref{s!lt3} (dimension~7, codimension~$d+2$) and $W$ in
\S\ref{s!k45} (dimension~10, codimension~6) present new Gorenstein rings
purpose built to model certain 3-fold flips of Type~A.

\section{The apolar variety $V_k$}\label{s!1}

The apolar varieties $V_k\subset\Aff^{k+5}$ introduced here provide  an
infinite family of affine Gorenstein 5-folds that are almost homogeneous
spaces under $\GL(2)\times\G_m$. We treat the $V_k$ as varieties in
their own right from several different points of view.

\subsection{The definition by equations}

We define 5-folds $V_k\subset\Aff^{k+5}_{\Span{x_{0\dots k},a,b,c,z}}$
for each $k\ge3$. First set up $2\times k$ and $k\times(k-2)$ matrixes
\begin{equation*}
M=
\begin{pmatrix}
x_0 & \dots & x_{i-1} & \dots & x_{k-1} \\
x_1 & \dots & x_i & \dots & x_k \\
\end{pmatrix}
%\quad\hbox{and}\quad
\end{equation*}
and
\begin{equation*}
N = \begin{pmatrix}
a & \\
b & a \\
c & b & a \\
& \vdots &&&& \vdots\\
&&&& c & b & a \\
&&&& & c & b \\
&&&&&& c
\end{pmatrix}.
\end{equation*}
Our variety $V_k\subset\Aff^{k+5}_{\Span{x_{0\dots k},z,a,b,c}}$ is
defined by two sets of equations:
\begin{equation}
\mathrm{(I)}\quad MN=0
\qquad\hbox{and \qquad (II)} \quad
\bigwedge^2M=z\cdot\bigwedge^{k-2}N.
\label{eq!Vk}
\end{equation}
(I) is a recurrence relation
\begin{equation}
ax_{i-1}+bx_i+cx_{i+1}=0 \quad\hbox{for $i=1,\dots,k-1$.}
\label{eq!rr}
\end{equation}
(II) is a $(k-2)\times k$ adaptation of Cramer's rule giving the
Pl\"ucker coordinates of the space of solutions of (I) up to a scalar
factor $z$. The order and signs of the minors in (II) is not a problem
here, as one sees from the guiding cases
\begin{equation*}
x_{i-1}x_{i+1}-x_i^2 = a^{i-1}c^{k-i-1}z \enspace\hbox{and}\enspace
x_{i-1}x_{i+2}-x_ix_{i+1} = a^{i-1}bc^{k-i-2}z.
%\label{eq!min}
\end{equation*}
(However, in subsequent cases, in particular when we work with Pfaffians
in \ref{ss!pf}, we need to fix a convention on their order and signs.)
Note that the maximal $(k-2)\times(k-2)$ minors of $N$ include $a^{k-2}$
(delete the last two row) and $c^{k-2}$ (delete the first two). More
generally, deleting two adjacent rows $i-1,i$ gives $a^{i-1}c^{k-i-1}$
as a minor (only the diagonal contributes), whereas deleting two rows
$i-1,i+1$ gives the minor $a^{i-1}bc^{k-i-2}$.

Thus our second set of equations is
\begin{equation*}
x_{i-1}x_{j+1}-x_ix_j = z\det N(i-1,j).
\end{equation*}
Relations for $x_ix_j-x_kx_l$ for all $i+j=k+l$ are obtained as
combinations of these; for example
\begin{equation*}
\begin{aligned}
x_{i-1}x_{j+2}-x_{i+1}x_j &=
x_{i-1}x_{j+2}-x_ix_{j+1} + x_ix_{j+1}-x_{i+1}x_j \\
& = zN(i-1,j+1)+zN(i,j).
\end{aligned}
\end{equation*}

\begin{theorem}
\label{th!Vk}
For $k\ge3$, (I) and (II) define a reduced irreducible
Gorenstein $5$-fold
\begin{equation*}
V_k\subset\Aff^{k+5}_{\Span{x_{0\dots k},a,b,c,z}}.
\end{equation*}
This also holds for $k=2$, with (II) involving interpreting the
$0\times0$ minors as the single equation $1\cdot z=x_0x_2-x_1^2$.
\end{theorem}

This theorem follows at once from the following lemma.

\begin{lem}
\label{lem!Vk}
\begin{enumerate}
\renewcommand{\labelenumi}{(\roman{enumi})}
\item $z$ is a regular element for $V_k$.

\item The section $z=0$ of $V_k$ is the quotient of the hypersurface
\begin{equation*}
\widetilde W:(g:=au^2+buv+cv^2=0) \subset \Aff^5_{\Span{a,b,c,u,v}}
\end{equation*}
by the $\bmu_k$ action $\frac1k(0,0,0,1,1)$. It is Gorenstein because
\begin{equation*}
\frac{\dd a \wedge\dd b \wedge\dd c \wedge\dd u \wedge\dd v}{g} \in
\om_{\Aff^5}(\widetilde{W}).
\end{equation*}
is $\bmu_k$ invariant.

\item Also $z,a,c$ is a regular sequence, and the section $z=a=c=0$ of
$V_k$ is the toric Gorenstein surface (three-sided tent) consisting of
$\frac1k(1,1)$ with coordinates $x_0,\dots,x_k$ and two copies of
$\Aff^2$ with coordinates $x_0,b$ and $x_k,b$.

\end{enumerate}
\end{lem}

\begin{pf} First, if $c\ne0$ then $a,b,c,x_0,x_1$ are free parameters,
and the recurrence relation (I) gives $x_2,\dots,x_k$ as rational
function of these. One checks that the first equation in (II) gives
$z=-\frac{ax_0^2+bx_0x_1+cx_1^2}{c^{k-1}}$ and the remainder follow.
Similarly if $a\ne0$.

If $a=c=0$ and $b\ne0$ then one checks that $x_0,x_k,b$ are free
parameters, $x_i=0$ for $i=1,\dots,k-1$ and $z=\frac{x_0x_k}{b^{k-2}}$.
Finally, if $a=b=c=0$ then $x_0,\dots,x_k$ and $z$ obviously parametrise
$\frac1k(1,1)\times\Aff^1$.

Therefore, no component of $V_k$ is contained in $z=0$, proving~(i).

After we set $z=0$, the equations (II) become $\bigwedge^2M=0$, and
define the cyclic quotient singularity $\frac1k(1,1)$ (the cone over the
rational normal curve). Introducing $u,v$ as the roots of
$x_0,\dots,x_k$, with $x_i=u^{k-i}v^i$, boils the equations $MN=0$ down
to the single equation $g:=au^2+buv+cv^2=0$. This proves (ii). (iii) is
easy.
\end{pf}

\subsection{The equations as Pfaffians} \label{ss!pf}

The equations of $V_k$ fit together as $4\times4$ Pfaffians of a skew
matrix. For this, edit $M$ and $N$ to get two new matrixes,
\begin{equation}
\renewcommand{\arraystretch}{1.2}
M'=
\begin{pmatrix}
x_0 & \dots & x_{i-1} & x_i && \dots & x_{k-2} \\
x_1 & \dots & x_i & x_{i+1} && \dots & x_{k-1} \\
x_2 & \dots & x_{i+1} & x_{i+2} && \dots & x_k
\end{pmatrix}
\label{eq!Mpr}
\end{equation}
which is $3\times(k-1)$ and $N'$, the $(k-1)\times(k-3)$ matrix with the
same display as $N$ (that is, delete the first (or last) row and column
of $N$). Equations (I) can be rewritten $(a,b,c)M'=0$.

Now all of the equations (\ref{eq!Vk}) can be written as the $4\times4$
Pfaffians of the $(k+2)\times(k+2)$ skew matrix
\begin{equation}
\renewcommand{\arraystretch}{1.0}
\renewcommand{\arraycolsep}{0.25em}
\begin{pmatrix}
\begin{matrix} c\; & -b \\
 & a \end{matrix}
 & M' \\[12pt]
& z\bigwedge^{k-3}N'
\end{pmatrix}.
\label{eq!kplus2skew}
\end{equation}
The Pfaffians $\Pf_{12.3(i+3)}$ give the recurrence relation
\eqref{eq!rr}, while the remaining Pfaffians give (II). In more detail,
the big matrix is
\begin{equation*}
\renewcommand{\arraystretch}{1.5}
\begin{pmatrix}
c & -b & x_0 & \dots & x_{i-1} & x_i && \dots & x_{k-2} \\
 & a & x_1 & \dots & x_i & x_{i+1} && \dots & x_{k-1} \\
  && x_2 & \dots & x_{i+1} & x_{i+2} && \dots & x_k \\
   &&& zc^{k-3} & \dots &&& \dots & \dots \\
    &&&& zc^{k-i-1}a^{i-2} & -zbc^{k-i-2}a^{i-2} && \dots & \dots \\
     &&&&& zc^{k-i-2}a^{i-1} && \dots & \dots \\
      &&&&&& & \dots & \dots \\
       &&&&&&&& za^{k-3}
\end{pmatrix}
\end{equation*}
with bottom right $(k-1)\times(k-1)$ block equal to the $(k-3)$rd wedge
of $N'$ (with signs).

\medskip\paragraph{\bf Small values of $k$.} Our family starts with $k\ge3$; the case
$k=2$ would give the hypersurface $ax_0+bx_1+cx_2=0$, with
$z:=x_0x_2-x_1^2$. The first regular case is $k=3$, which gives the
$5\times5$ skew determinantal
\begin{equation*}
\renewcommand{\arraystretch}{1.3}
\begin{pmatrix}
c & -b & x_0 & x_1 \\
 & a & x_1 & x_2 \\
 && x_2 & x_3 \\
 &&& z
\end{pmatrix}
\end{equation*}
a regular section of the affine Grassmannian $\aGr(2,5)$. The case $k=4$
is
\begin{equation*}
\renewcommand{\arraystretch}{1.3}
\begin{pmatrix}
c & -b & x_0 & x_1 & x_2 \\
 & a & x_1 & x_2 & x_3 \\
 && x_2 & x_3 & x_4 \\
 &&& zc & -zb \\
 &&&& za
\end{pmatrix},
\end{equation*}
an easy case of the standard extrasymmetric $6\times6$ determinantal of
Dicks and Reid, \cite{TJ}, 9.1, equation (9.4).

The first really new case is $k=5$, with equations the $4\times4$
Pfaffians of the $7\times7$ skew matrix
\begin{equation*}
\renewcommand{\arraystretch}{1.3}
\renewcommand{\arraycolsep}{0.4em}
\begin{pmatrix}
c & -b & x_0 & x_1 & x_2 & x_3 \\
 & a & x_1 & x_2 & x_3 & x_4 \\
 && x_2 & x_3 & x_4 & x_5 \\
 &&& zc^2 & -zbc & z(b^2-ac) \\
 &&&& zac & -zab \\
 &&&&& za^2
\end{pmatrix}
%\label{eq!110}
\end{equation*}
We first arrived at this matrix by guesswork (with the $z$ floated over
from the row-columns $4,5,6,7$ to $1,2,3$), determining the
superdiagonal entries $c^2,ac,a^2$ and those immediately above $-bc,-ac$
by eliminating variables to smaller cases; the entry $b^2-ac$ is then
fixed so that the bottom $4\times4$ Pfaffian vanishes identically.

\subsection{The variety $V_k$ by apolarity}\label{ss!Vk}
We can treat $V_k$ as an almost homogeneous space under
$\GL(2)\times\G_m$. For this, view $x_0,\dots,x_k$ as coefficients of a
binary form and $a,b,c$ as coefficients of a binary quadratic form in
dual variables, so that the equations $MN=0$ or $(a,b,c)M'=0$ are the
apolarity relations. In general terms, {\em polarity} can be described
as a choice of splitting of maps such as $\Sym^{d-1}U \tensor U\onto
\Sym^dU$ (here $U=\C^2$ is the given representation of $\GL(2)$), or more
vaguely as a way of viewing the $2\times d$ matrix
$\bigl(\begin{smallmatrix} y_0 & \dots & y_{d-1} \\ y_1 & \dots &
y_d\end{smallmatrix}\bigr)$ or his bigger cousin \eqref{eq!Mpr} as a
single object in determinantal constructions.

More formally, write
\begin{equation*}
\begin{aligned}
q &= a\checku^2+2b\checku \checkv+c\checkv^2 \in \Sym^2U^\vee \quad\hbox{and} \\
f &= x_0u^k+kx_1u^{k-1}v+{\textstyle{\binom k 2}}x_2 u^{k-2}v^2+\dots+x_kv^k \in \Sym^kU.
\end{aligned}
\end{equation*}
Including the factor $\binom k i$ in the coefficient of $u^iv^{k-i}$ is
a standard move in this game.

The second polar of $f$ is the polynomial
\begin{align*}
\Phi(u,v,u',v')&=\frac1{k(k-1)}\left(\frac{\partial^2 f}{\partial u^2}\tensor u'{}^2
+2\frac{\partial^2 f}{\partial u \partial v}\tensor u'v'
+\frac{\partial^2 f}{\partial v^2}\tensor v'{}^2\right)
\\
&=
\sum_{i=0}^{k-2} \textstyle\binom {k-2}ix_iu^{k-i-2}v^i\tensor u'{}^2\\[-5pt]
&\kern1.6cm+2\sum_{i=1}^{k-1} \textstyle\binom {k-2}{i-1} x_iu^{k-i-1}v^{i-1}\tensor u'v' \\[-5pt]
&\kern3.2cm+\sum_{i=2}^{k} \textstyle\binom {k-2}{i-2}x_iu^{k-i}v^{i-2}\tensor v'{}^2
\\
&=
\sum_{i=0}^{k-2} \textstyle\binom {k-2}iu^{k-2-i}v^i\tensor
 \bigl(x_iu'{}^2 + 2x_{i+1}u'v' + x_{i+2}v'{}^2 \bigr)\\
&\kern1.6cm \in \Sym^{d-2}U\tensor\Sym^2U.
\end{align*}
We apply $q\in\Sym^2U^\vee$ to the second factor and equate to zero to
obtain the recurrence relation $(a,b,c)M=0$. In other words, substitute
$u'{}^2\mapsto a$, $u'v'\mapsto \frac12 b$, and $v'{}^2\mapsto c$ in
$\Phi$.

Moreover, the second set of equations follows from the first by
substitution, provided (say) that $c\ne0$ and we fix the value of
$x_0x_2-x_1^2$; for example, in
\begin{equation*}
x_ix_{i+2}-x_{i+1}^2
\end{equation*}
substituting $x_{i+2}=-\frac ac x_i-\frac bc x_{i+1}$ gives
\begin{equation*}
x_i\Bigl(-\frac ac x_i-\frac bc x_{i+1}\Bigr)-x_{i+1}^2
 =-\frac ac x_i^2-\Bigl(\frac bc x_i+x_{i+1}\Bigr)x_{i+1},
\end{equation*}
and we can substitute $-\frac ac x_{i-1}$ for the bracketed expression,
to deduce that
\begin{equation*}
x_ix_{i+2}-x_{i+1}^2 = \frac ac\left(x_{i-1}x_{i+1}-x_i^2\right), \quad\hbox{etc.}
\end{equation*}

A normal form for a quadratic form under $\GL(2)$ is $uv$, so that a
typical solution to the equations is
\begin{equation*}
(a,b,c) = (0,1,0), \quad (x_{0\dots k})=(1,0,\dots,0,1), \quad z=1.
\end{equation*}
in the representation $\Sym^2U^\vee \oplus \Sym^kU \oplus \C^1$ of
$\GL(2)\times\G_m$, where the final $\G_m$ acts by homotheties on
$U^\vee$, so acts on $q\in \Sym^2U^\vee$ by $q\mapsto\la^2q$ and on $z$
by $z\mapsto\la^2z$. Then $V_k$ is the closure of the orbit of this
typical apolar vector.

\section{Diptych varieties and Mori flips of Type A}
\label{s!dip1}

The varieties $V_k\subset\Aff^{k+5}$ form a simple and natural series
of Gorenstein 5-folds, each with an action of a large algebraic group and,
by Lemma~\ref{lem!Vk}, a regular sequence $z,a,c$ whose common zero locus
is a reducible toric surface composed of a cycle of three affine toric surfaces.

In \cite{BR1}, we introduce a rather more complicated series of Gorenstein varieties:
these are 6-folds
\[
V_{ABLM}\subset\Aff^{k+l+6}
\]
(where $l$ is the number appearing in \eqref{eq!defn_of_l}), each
admitting a regular sequence $A,B,L,M$ whose common zero locus $T\subset
V_{ABLM}$ is a reducible toric surface composed of a cycle of four
affine toric surfaces which we call a {\em tent}. There is more
combinatorial structure inside $V_{ABLM}$: namely $V_{LM}:=(A=B=0)$ and
$V_{AB}:=(L=M=0)$ are toric 4-folds inside $V_{ABLM}$ whose intersection
equals $T$. In the language of \cite{AH}, $V_{ABLM}$ is an affine
T-variety (T for torus, not for tent): it admits an action of a torus
$\T=(\Gm^\times)^4$ which restricts to the intrinsic torus action on
each of the toric strata described so far.

Each diptych variety depends on a 2-step recurrent continued fraction
$[d,e,d,\dots]$ to $k$ terms. Starting from nothing, this data
determines the toric configuration $V_{AB}\supset T \subset V_{LM}$, and
the {\em existence of diptych varieties} is then the claim that this
configuration arises inside an irreducible 6-fold, the diptych variety,
as above; this claim is proved in the case $de>4$, $d,e\ge2$ in
\cite{BR1}.

In \S\ref{s!dip} we use $V_k$ to prove the existence of diptych
varieties in the case $de=4$.
We need some of the definitions and notions of \cite{BR1} for this.
Given integers $d,e,k\ge1$, consider the continued fraction expansion
with $k$ terms
\[
[d,e,d,\dots] = d - \frac{1}{e-\cdots}.
\]
Define $[b_1,\dots,b_{l-1}]$ to be the complementary continued
fraction of a truncation as follows.
Truncate the expansion $[d,e,d,\dots]$ to $k-1$ terms and reverse it,
and then consider the uniquely defined minimal sequence of $b_j\ge2$
for which
\begin{equation}
[\dots,d,e,d,1,b_{l-1},\dots,b_1] = 0.
\label{eq!defn_of_l}
\end{equation}
For example, starting with $[4,3,4]$, one calculates
$[3,4,1,2,2,3,2]=0$, so in this case one has $[b_1,b_2,b_3,b_4] = [2,3,2,2]$.
(This is the Riemenschneider complementary continued fraction,
in the sense of \cite{BR1} Proposition~2.1(d).)
Set $b_l=1$.

Now define a toric variety $V_{AB}$ as follows.
Start with four variables $x_k,y_l$, $A,B$. Define the Laurent monomial
$x_{k-1} = Ax_k^dy_l^{-1}$, and then
\begin{equation}
\label{eq!xtags}
x_{k-2i} = x_{k-2i+1}^e x_{k-2i+2}^{-1}
\quad\hbox{and}\quad
x_{k-2i-1} = x_{k-2i}^d x_{k-2i+1}^{-1}
\end{equation}
alternating the exponents $d,e$ until you reach $x_0$.
Similarly define $y_{l-1}=Bx_k^{-1}y_l^{b_l}$, and then
\[
y_{j-1} = y_j^{b_j}y_{j+1}^{-1}
\]
until you reach $y_0$.
We treat these expressions in two ways: first as monomials
in a lattice $\M_{AB} = \Z^4$ based by $A,B,x_k,y_l$;
second as independent variables $A,B$,
$x_{0\dots k}$, $y_{0\dots l}$ on affine space $\Aff^{k+l+4}$.
The cone
\[
\sigma_{AB} = \left< A,B,x_0,\dots, x_k, y_0,\dots, y_l \right>
\subset \M_{AB}
\]
defines a toric variety $V_{AB} = X_{\sigma_{AB}}$ which embeds
naturally as
\[
V_{AB} \subset \Aff^{k+l+4}
\]
defined by the relations above (after multiplying up denominators)
and others that follow from syzygies. (In other words, the relations
above define a union of components, of which $V_{AB}$ is the
unique component not contained in a coordinate hyperplane.)
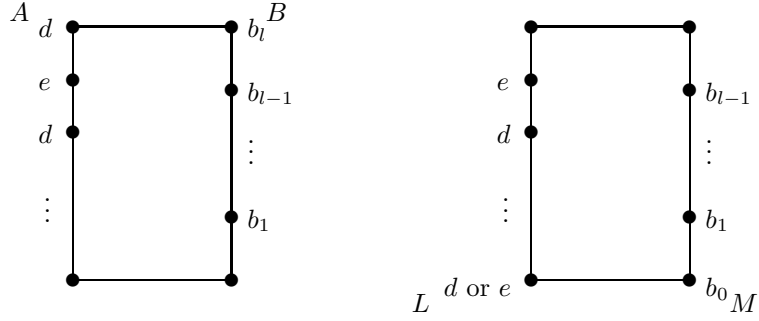
\begin{figure}[ht]
\begin{picture}(0,110)(-50,0) % 80 -> -50 for page width
\put(45,7){\line(0,1){96}}
\put(105,7){\line(0,1){96}}
\put(45,7){\line(1,0){60}}
\put(45,103){\line(1,0){60}}
\put(45,7){\circle*5}
\put(45,63){\circle*5}
\put(45,83){\circle*5}
\put(45,103){\circle*5}
%% \put(15,-4){$(-(b_1-1))$}
\put(34,29){$\vdots$}
\put(32,59){$d$}
\put(32,79){$e$}
\put(32,99){$d$}
\put(21,105){$A$}
\put(105,7){\circle*5}
\put(105,31){\circle*5}
%\put(105,55){\circle*5}
\put(105,79){\circle*5}
\put(105,103){\circle*5}
%% \put(111,1){(0)}
\put(111,27){$b_1$}
\put(111,51){$\vdots$}
\put(111,75){$b_{l-1}$}
\put(111,99){$b_l$}
\put(118,105){$B$}
\end{picture}
\begin{picture}(150,110)(-220,0) % -80 -> -220 for page width
\put(45,7){\line(0,1){96}}
\put(105,7){\line(0,1){96}}
\put(45,7){\line(1,0){60}}
\put(45,103){\line(1,0){60}}
\put(45,7){\circle*5}
\put(45,63){\circle*5}
\put(45,83){\circle*5}
\put(45,103){\circle*5}
\put(12,1){$d$ or $e$}
\put(34,29){$\vdots$}
\put(32,59){$d$}
\put(32,79){$e$}
%% \put(-17,99){$(-(b_{l-1}-1))$}
\put(0,-5){$L$}
\put(105,7){\circle*5}
\put(105,31){\circle*5}
%\put(105,55){\circle*5}
\put(105,79){\circle*5}
\put(105,103){\circle*5}
\put(111,1){$b_0$}
\put(111,27){$b_1$}
\put(111,51){$\vdots$}
\put(111,75){$b_{l-1}$}
%% \put(111,99){(0)}
\put(120,-5){$M$}
\end{picture}
\caption{The pair of long rectangles for $[d,e,d,\dots]$ to $k$ terms}
\label{f!21}
\end{figure}

Similarly we define $V_{LM}$ starting from the four variables $x_0,y_0$,
$L,M$ and applying analogous relations for $x_1,x_2,\dots$ and
$y_1,y_2,\dots$ but with the terms of the reversed continued fraction:
that is, with $[d,e,d,\dots]$
if $k$ is even, and from $[e,d,e,\dots]$ to $k$ terms if $k$ is odd.
Again there is a lattice $\M_{LM}$ containing the defining cone $\sigma_{LM}$.

We sketch all of this data in a picture, called a {\em pair of long rectangles},
as in Figure~\ref{f!21}, in which the bullet points represent $x_0,x_1,\dots,x_k$
up the left-hand side of each long rectangle and $y_0,\dots,y_l$ up the
right-hand side, the {\em tags} $d$, $e$ and $b_j$ appear next to the
corresponding variable on which they appear as an exponent,
and the four auxilliary variables, or {\em annotations}, $A,B,L,M$
positioned near the corners where they appear in the initial defining relations.
Influenced by this picture, we refer to data associated to $x_0, y_0$ as
the {\em bottom end} of the long rectangles, and to $x_k, y_l$ as
the {\em top end}.

Notice from the defining relations that the lattices $\M_{AB}$ and
$\M_{LM}$ are in fact identical, and so we identify them as $\M$. To
avoid prejudice, we use the {\em impartial basis} $L,M,A,B$ of $\M$.
Although these four monomials are only a $\Q$-basis spanning an index
$de$ sublattice of $\M$, expressing lattice points in them turn out to
express the antagonistic convexity properties of $\sigma_{AB}$ and
$\sigma_{LM}$ most cleanly.

Although it is not completely obvious, the data assembled so far
describes the toric monomial cones of the configuration $V_{AB}\supset T
\subset V_{LM}$ for the initial continued fraction expansion
$[d,e,d,\dots]$; see \cite{BR1}, \S3. To show the existence of the
corresponding diptych \hbox{6-fold}, we simply build its equations from
the bottom end up. We start by combining the equations of $V_{AB}$ and
$V_{LM}$ at the bottom end in a naive way:
\begin{align}
x_1y_0 &= y_1A^\alpha B^\beta + x_0^{(d\hbox{ or }e)}L \nonumber \\
x_0y_1 &= A^\gamma B^\delta + y_0M, \label{eq!x0y0}
\end{align}
where the exponents $\alpha$, $\beta$, $\lambda$, $\mu$ are determined
by the tag relations we started from (and, unsurprisingly, appear in
convergents of the continued fraction expansion $[d,e,d,\dots]$). These
relations define a Gorenstein \hbox{6-fold}
$V_0\subset\Aff^8_{\Span{A,B,L,M,x_0,x_1,y_0,y_1}}$, that contains a
divisor
\[
D_0 = (x_0=y_0 = A^\lambda B^\mu = 0) \subset V_0,
\]
where $A^\lambda B^\mu = \gcd(A^\alpha B^\beta, A^\gamma B^\delta)$. We
now apply the Gorenstein unprojection theorem of \cite{PR} serially to
construct a sequence of pairs $D_\nu\subset V_\nu$, adding the remaining
variables $x_i$, $y_j$ one at a time until we reach $V_\nu=V_{ABLM}$.

We demonstrate the first step by use of a magic {\em pentagram}:
we seek to include the variable $x_2$ and calculate any relations that involve it.
Consider the $5\times 5$ antisymmetric matrix (we write only the strict
upper triangle), which we also refer to as the {\em Pfaffian matrix},
\begin{equation}
\label{eq!pfaff}
M_0=
\begin{pmatrix}
x_1 & y_1A^{\alpha-\lambda} B^{\beta-\mu} & -x_0^{(d\hbox{ or }e)-1}L & -x_2 \\
& x_0 & A^\lambda B^\mu & -M\\
&& y_0 & A^{\gamma-\lambda} B^{\delta-\mu} \\
&&& y_1
\end{pmatrix}.
\end{equation}
The first and last of the maximal Pfaffians of $M$ give precisely the pair of
relations~\eqref{eq!x0y0}.
The other three maximal Pfaffians involve expressions for $x_2\cdot I_{D_0}$,
where $I_{D_0}=(x_0,y_0,A^\lambda B^\mu)$ is the defining
ideal of the unprojection divisor $D_0\subset V_0$.
These five Pfaffians define a Gorenstein variety $V_1\subset\Aff^9$
in variables $A,B,L,M$, $x_0,x_1,x_2,y_0,y_1$. If $k=1$, then
this is $V_{ABLM}$, otherwise it contains a divisor
\[
D_1=(x_0=x_1=y_0=A^? B^? = 0) \subset V_1,
\]
where the exponents on $A^?B^?$ can be determined from the particular
values of $d,e,k$. One can check that the 4-fold locuses $(A=B=0)$ and
$(L=M=0)$ and their surface intersection correspond to the toric
configuration; this is part of the claim of the existence of diptych
varieties. The five equations constructed here have leading terms
\begin{gather*}
x_0y_1 = \cdots \qquad
x_1y_0 = \cdots \\
x_2x_0 = \cdots \qquad
x_2y_0 = \cdots \qquad
x_1y_1 = \cdots,
\end{gather*}
and joining these pairs of variables on Figure~\ref{f!21}
draws a pentagram -- hence the name. (It is magic because
it works.)

The order we add the variables is important. We lay a {\em bar} at the
level of variables we have considered so far: we start with the bar
$x_1\frac{\qquad}{}y_1$, to indicate that we have all variables below
these, then raise it to $x_2\frac{\qquad}{}y_1$ and so on as we add
subsequent variables. Fortunately the precise order required is a
technical point that our use of $V_k$ in this paper sidesteps.

As an exercise, one can write an alternative proof of Theorem~\ref{th!Vk} above
in the style of \cite{BR1}:
start with any of the codimension 2 complete intersections
\begin{equation*}
\begin{pmatrix}
x_{i-1}x_{i+1}=x_i^2+a^{i-1}c^{k-i-1}z \\
ax_{i-1}+bx_i+cx_{i+1}=0
\end{pmatrix}
\subset
\Aff^7_{\Span{x_{i-1},x_i,x_{i+1},a,b,c,z}}
\end{equation*}
and add the remaining variables one at a time as a serial unprojection
using magic pentagrams at each step. (Or see \cite{BR1}, 1.2, for
a fully-worked example of a similar calculation.)

Once set up properly, much of this construction is automatic. Curiously,
the hardest part, and the bulk of the subtle machinery developed in
\cite{BR1}, is to show that the natural unprojection divisor $D_\nu$ is
a subscheme of $V_\nu$. Again, our use of the $V_k$ here completely
sidesteps that point -- when we need to make unprojection arguments in
\S\ref{s!dip}, the inclusion of the divisor is straightforward.

The contrast between the simple geometric constructions of this paper
and the delicate and lengthy methods of \cite{BR1} is striking. The
varieties $V_k$ arise naturally from the representation theory of
$\GL(2)\times\G_m$, in contrast to any construction we could find in
\cite{BR1}. There is still some work to do in Section~\ref{s!dip} to go
from $V_k$ to the diptych varieties, but it is easy compared to
\cite{BR1}. Whether the other dip\-tychs of \cite{BR1} can be modelled
on almost homogeneous spaces in a similar way remains a mystery; this
point has eluded us for a couple of decades.

\section{Application of $V_k$ to diptych varieties with $de=4$}
\label{s!dip}
Diptych varieties $V_{ABLM}$ depend on three parameters $d,e,k\ge1$. The
solutions of $de=4$ are $(d,e)=(2,2)$, $(4,1)$ and $(1,4)$, and we allow
any $k\ge1$. In each case, we construct almost all of the coordinate
ring of $V_{ABLM}$ by a regular pullback from the key variety $V_k$ of
\S\ref{s!1}. We then adjoin the remaining few variables by an
unprojection argument using the ideas of \S\ref{s!dip1}. Our proofs here
are selfcontained, but we refer to \cite{BR1} in places this clarifies
the argument; see especially the worked example \cite{BR1},~1.2.)

\subsection{Case $[2,2]$}
We first construct the diptych variety $V_{ABLM}$ with the monomial
cones $\si_{AB}$ and $\si_{LM}$ of Figure~\ref{f!22}.
\begin{figure}[ht]
\begin{picture}(0,110)(-50,0)
\put(45,7){\line(0,1){96}}
\put(105,7){\line(0,1){96}}
\put(45,7){\line(1,0){60}}
\put(45,103){\line(1,0){60}}
\put(45,7){\circle*5}
\put(45,23){\circle*5}
\put(45,83){\circle*5}
\put(45,103){\circle*5}
\put(24,1){(0)}
\put(34,49){$\vdots$}
\put(32,19){2}
\put(32,79){2}
\put(32,99){2}
\put(21,105){$A$}
\put(105,7){\circle*5}
% \put(105,31){\circle*5}
\put(105,55){\circle*5}
% \put(105,79){\circle*5}
\put(105,103){\circle*5}
\put(111,1){$(-1)$}
% \put(111,27){2}
\put(111,51){$k$}
% \put(111,75){2}
\put(111,99){1}
\put(118,105){$B$}
\end{picture}
\begin{picture}(150,110)(-220,0)
\put(45,7){\line(0,1){96}}
\put(105,7){\line(0,1){96}}
\put(45,7){\line(1,0){60}}
\put(45,103){\line(1,0){60}}
\put(45,7){\circle*5}
\put(45,23){\circle*5}
\put(45,83){\circle*5}
\put(45,103){\circle*5}
\put(32,1){2}
\put(34,49){$\vdots$}
\put(32,19){2}
\put(32,79){2}
\put(24,99){(0)}
\put(21,-5){$L$}
\put(105,7){\circle*5}
% \put(105,31){\circle*5}
\put(105,55){\circle*5}
% \put(105,79){\circle*5}
\put(105,103){\circle*5}
\put(111,1){1}
% \put(111,27){2}
\put(111,51){$k$}
% \put(111,75){2}
\put(111,99){$(-1)$}
\put(118,-5){$M$}
\end{picture}
\caption{Case $[2,2]$}
\label{f!22}
\end{figure}
It has variables $x_{0\dots k}$ on the left against $y_{0\dots 2}$ on
the right, tagged as in Figure~\ref{f!22}, together with $A,B$, $L,M$.
Although we do not yet own $V_{ABLM}$, we know some of its equations:
by \eqref{eq!x0y0}, we find the two bottom equations:
\begin{equation}
x_1y_0= A^{k-1}B^k+x_0^2L
\quad\hbox{and}\quad
x_0y_1= ABx_1+y_0M.
\label{eq!V0}
\end{equation}
Then, following the model of \eqref{eq!pfaff},
the pentagram $y_1,y_0,x_0,x_1,x_2$ adjoins $x_2$, and
$x_3,\dots,x_k$ are adjoined by a long rally of {\em flat} pentagrams
$y_1,x_{i-1},x_i,x_{i+1},x_{i+2}$ with matrixes
\begin{equation}
\renewcommand{\arraystretch}{1.3}
\renewcommand{\arraycolsep}{0.3em}
\begin{pmatrix}
y_1 & x_1 & -M & -x_2 \\
& y_0 & AB & -x_0L \\
&& x_0 & A^{k-2}B^{k-1} \\
&&& x_1
\end{pmatrix}
\label{eq!pent}
\end{equation}
and
\begin{equation*}
%\enspace \hbox{and} \enspace
\begin{pmatrix}
y_1 & x_{i+1} & -LM & -x_{i+2} \\
& x_{i-1} & AB & -x_i \\
&& x_i & (AB)^{k-i-2}(LM)^{i-1}BM \\
&&& x_{i+1}
\end{pmatrix}
\end{equation*}
giving the Pfaffian equations
\begin{gather*}
y_1x_i=ABx_{i+1}+LMx_{i-1}, \quad
x_{i-1}x_{i+1}=x_i^2+(AB)^{k-i-1}(LM)^{i-1}BM \\[4pt]
\hbox{and}\quad
x_{i-1}x_{i+2}=x_ix_{i+1}+(AB)^{k-i-2}(LM)^{i-1}BMy_1.
\end{gather*}
We see that these are the equations of $V_k$ after the substitution
\begin{equation}
\label{eq!subst}
(a,b,c,z)\mapsto(LM,-y_1,AB,BM).
\end{equation}
Thus to construct our diptych variety, we pull back
$V_k\subset\Aff^{k+5}$ by \eqref{eq!subst}, then adjoin the two corners
$y_0$, $y_2$ as unprojection variables. Adjoining either of these is
easy, but adjoining the second then requires a simple application of
some of the main ideas of proof in Sections~4--5 of \cite{BR1} which
we work out here.

\begin{lem}
\label{lem!W1}
Define $W_0\subset \Aff^{k+6}_{\Span{x_{0\dots k},y_1,A,B,L,M}}$ to be
the pullback of\/ $V_k$ under the morphism $\Aff^{k+6}\to\Aff^{k+5}$ given
by
\eqref{eq!subst}.
\begin{enumerate}
\renewcommand{\labelenumi}{(\roman{enumi})}
\item $W_0\subset \Aff^{k+6}$ is an irreducible $6$-fold.

\item $D_0=(x_1=\cdots=x_k=M=0)$ is contained in $W_0$ as a divisor.

\item The unprojection $W_1\subset \Aff^{k+6}\times\Aff^1_{\Span{y_0}}$
of $D_0\subset W_0$ with unprojection variable $y_0$ includes the
equations \eqref{eq!V0} as generators of its defining ideal.
\end{enumerate}
\end{lem}

\begin{pf} (ii) is immediate from the defining equations \eqref{eq!Vk}
of $V_k$: setting $x_1=\cdots=x_k=0$ leaves only terms divisible by
$M$. It is a divisor because it has the right dimension.
(iii) follows from the Pfaffians of the matrix
\eqref{eq!pent}, that express the unprojection variable $y_0$ as a
rational function in $x_0,x_1,y_1$, $A,B,L,M$ with a simple pole on $D$.
This includes the equations \eqref{eq!V0}.
\end{pf}

Once we own $y_0\in \C[W_1]$, we have to establish that the unprojection
divisor of $y_2$ is contained in the variety $W_1$. The detailed
statement is Theorem~\ref{th!22} below. (This is the same as the key point of
the proof of \cite{BR1}, but our case here is much easier.) To prove it,
we work with the $\T$-weights of each homogeneous polynomial in
$x_0,\dots,y_2$, $A,B,L,M$, written in terms of the impartial basis dual to the monomials $L,M,A,B$
(compare \cite{BR1}, Proposition~4.1).
These base a slightly smaller lattice, giving some of
the impartial coordinates of monomials little denominators $d$ or $e$.
The tag equations of $V_{AB}$ and $V_{LM}$ from Figure~\ref{f!22}
determine the impartial coordinates, as follows.

\begin{lem}
\label{lem!pp}
In the impartial basis $L,M,A,B$, the monomials $x_0$, \dots, $y_2$
have $\T$-weights:
\begin{equation*}
\renewcommand{\arraystretch}{1.2}
\renewcommand{\arraycolsep}{.2em}
\begin{array}{rclccccc}
&&&L&M&A&B\\[4pt]
x_0 &=& (&-\frac12&0&\frac{k-1}2&\frac k2&) \\
x_1 &=& (&0&\frac12&\frac{k-2}2&\frac{k-1}2&) \\
x_2 &=& (&\frac12&1&\frac{k-3}2&\frac{k-2}2&) \\
 &\vdots \\
x_i &=& (&\frac{i-1}2&\frac i2&\frac{k-i-1}2&\frac{k-i}2&) \\
 &\vdots \\
x_{k-1} &=& (&\frac{k-2}2&\frac{k-1}2&0&\frac12&) \\
x_k &=& (&\frac{k-1}2&\frac k2&-\frac12&0&)
\end{array}
\ \hbox{and}\quad
\renewcommand{\arraystretch}{1.2}
\renewcommand{\arraycolsep}{.35em}
\begin{array}{rccclccccc}
&&&L&M&A&B\\[4pt]
y_0 &=& (&0&-\frac12&\frac k2&\frac{k-1}2&)\\
y_1 &=& (&\frac12&\frac12& \frac12 & \frac12 &)\\
y_2 &=& (&\frac k2&\frac{k+1}2&0&-\frac12&) \\[20pt]
\end{array}
\end{equation*}
\end{lem}
\begin{pf}
These vectors satisfy all the tag relations of the pair of long
rectangles; or if you prefer, plug in the formulas from \cite{BR1},
Proposition~4.1.
\end{pf}

The following statement specifies the unprojection divisor $D_1\subset
W_1$ of $y_2$, completing our construction.

\begin{theorem}
\label{th!22}
In the notation of Lemma~\ref{lem!W1}, define
\begin{equation*}
D_1 = (x_0=\cdots=x_{k-1}=y_0=B=0)\subset
\Aff^{k+7}_{\Span{x_{0\dots k},y_0,y_1,A,B,L,M}}.
\end{equation*}
Then $D_1\subset W_1$, and the unprojection of $D_1$ in $W_1$ is the
diptych variety $V_{ABLM}$ on the pair of long rectangles of
Figure~\ref{f!22}.
\end{theorem}

\begin{pf} Most of the generators of $I_{W_1}$ are already in the ideal
of $I_{W_0}$, and so lie in the ideal $I_{D_1}$ by the argument of
Lemma~\ref{lem!W1} applied to $y_2$ rather than $y_0$. The equation
\eqref{eq!V0} of the form $x_0y_1=\cdots$ is known by
Lemma~\ref{lem!W1}(iii), and also lies in $I_{D_1}$.

The remaining generators of $I_{W_1}$ have leading terms $x_iy_0$ for
$i=1,\dots,k$. To prove that each of these lies in $I_{D_1}$, we prove a
stronger statement: every monomial in any of these generator
relations is divisible by one of $x_{0\dots k-1}$, $y_0$ or $B$.
In fact, we prove some stronger still. As in \cite{BR1},~5.1,
rather than working directly with these generators, we work with their
$\T$-weights, and we show that any monomial of $\T$-weight
equal to that of $x_iy_0$ (that is, any monomial that could appear
in a $\T$-homogeneous equation which included $x_iy_0$)
is divisible by one of $x_{0\dots k-1}$, $y_0$ or $B$.

For monomials $m,n$, write $m\simT n$ if $m$ and $n$ have the same
$\T$-weight, or equivalently, the same impartial coordinates. Suppose
$m\in \C[W_1]$ is a monomial with $m\simT x_iy_0$ for some $i=1,\dots,k$.
(Any term in the equation having leading term $x_iy_0$ satisfies this
equivalence, so if each such monomial lies in $I_{D_1}$ then certainly
the generator itself does.) We may assume that the monomial $m$ is of
the form $x_k^\xi y_2^\eta L^\la M^\mu A^\al B^\be$, since the other
variables already lie in $I_{D_1}$. We may assume further that $\xi=0$:
otherwise, dividing through by $x_i$, the
$\T$-weight of $y_0$ can be calculated from that of $(x_k/x_i)x_k^{\xi-1}$
times other variables whose $M$ coefficient is nonnegative; but this
has $M$ coefficient $>0$, whereas $y_0$ has $M$ coefficient $=-1/2$,
a contradiction.

Now compare $x_iy_0$ and $m=y_1^\eta L^\la M^\mu A^\al B^\be$: their
impartial coordinates are
\begin{equation*}
\renewcommand{\arraycolsep}{.4em}
\begin{array}{rcccccccc}
x_iy_0 &=& \bigl(\ & \frac{i-1}2&\frac{i-1}2&\frac{2k-i-1}2&\frac{2k-1+1}2 & \bigr)\phantom{.} \\[3pt]
y_1^\eta L^\la M^\mu A^\al B^\be
&=& \bigl(\ & \frac\eta2+\la&\frac\eta2+\mu& \frac\eta2+\al & \frac\eta2+\be & \bigr).
\end{array}
\end{equation*}
Since $\al\ge0$, it follows from the coefficient of $A$ that $\eta/2 \le (2k-i-1)/2$,
so now from the coefficient of $B$ we have $\be\ge1$.
In other words, $B$ divides the monomial $m$, and $m\in I_{D_1}$ as
required. \end{pf}

\subsection{Case $[4,1]$ with even $l=2k$}
\begin{figure}[ht]
\begin{picture}(0,110)(-50,0)
\put(45,7){\line(0,1){96}}
\put(105,7){\line(0,1){96}}
\put(45,7){\line(1,0){60}}
\put(45,103){\line(1,0){60}}
\put(45,7){\circle*5}
\put(45,63){\circle*5}
\put(45,83){\circle*5}
\put(45,103){\circle*5}
\put(15,1){$(-1)$}
\put(34,29){$\vdots$}
\put(32,59){2}
\put(32,79){2}
\put(32,99){3}
\put(21,105){$A$}
\put(105,7){\circle*5}
\put(105,31){\circle*5}
\put(105,55){\circle*5}
\put(105,79){\circle*5}
\put(105,103){\circle*5}
\put(111,1){(0)}
\put(111,27){2}
\put(111,51){$k+1$}
\put(111,75){2}
\put(111,99){1}
\put(118,105){$B$}
\end{picture}
\begin{picture}(150,110)(-220,0)
\put(45,7){\line(0,1){96}}
\put(105,7){\line(0,1){96}}
\put(45,7){\line(1,0){60}}
\put(45,103){\line(1,0){60}}
\put(45,7){\circle*5}
\put(45,63){\circle*5}
\put(45,83){\circle*5}
\put(45,103){\circle*5}
\put(32,1){3}
\put(34,29){$\vdots$}
\put(32,59){2}
\put(32,79){2}
\put(15,99){$(-1)$}
\put(21,-5){$L$}
\put(105,7){\circle*5}
\put(105,31){\circle*5}
\put(105,55){\circle*5}
\put(105,79){\circle*5}
\put(105,103){\circle*5}
\put(111,1){1}
\put(111,27){2}
\put(111,51){$k+1$}
\put(111,75){2}
\put(111,99){(0)}
\put(118,-5){$M$}
\end{picture}
\caption{Case $[4,1]$ with even $l=2k$}
\label{f!41}
\end{figure}
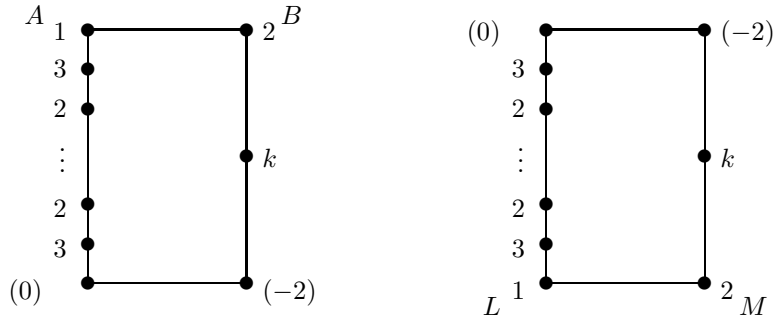

The odd numbered $x_i$ are redundant generators, and omitting them gives
Figure~\ref{f!41}. The diptych variety has variables $x_{0\dots k}$,
$y_{0\dots 4}$, $A,B$, $L,M$ with the two bottom equations
\begin{equation*}
x_1y_0= A^{k-1}B^{2k-1}y_1+x_0^3L
\quad\hbox{and}\quad
x_0y_1= A^kB^{2k+1}+y_0M.
\end{equation*}
We adjoin $y_2$, then $x_2,\dots,x_k$ by a game of pentagrams centred on
a long rally of flat pentagrams, with $y_2$ against
$x_{i-1},x_i,x_{i+1},x_{i+2}$ and Pfaffian equations
\begin{gather*}
y_2x_i=AB^2x_{i+1}+LM^2x_{i-1}, \\[4pt]
x_{i-1}x_{i+1}=x_i^2+(AB^2)^{k-i-1}(LM^2)^{i-1}BM \\[4pt]
\hbox{and}\quad
x_{i-1}x_{i+2}=x_ix_{i+1}+(AB^2)^{k-i-2}(LM^2)^{i-1}BMy_2
\end{gather*}
These are the equations of $V_k$ after the substitution
\begin{equation}
\label{eq!subst41}
(a,b,c,z)\mapsto(LM^2,-y_2,AB^2,BM).
\end{equation}

\begin{lem}
\label{lem!pp41}
In the impartial basis $L,M,A,B$, the monomials $x_0$, \dots, $y_4$ have
$\T$-weights as listed in Table~\ref{tab!LMAB}.
\end{lem}
%% \begin{equation*}
%% \renewcommand{\arraystretch}{1.2}
%% \renewcommand{\arraycolsep}{.2em}
%% \begin{array}{rclccccc}
%% x_0 &=& (&-\frac14&0& \frac{2k-1}4 & k &) \\
%% x_1 &=& (&\frac14&1& \frac{2k-3}4 & k-1 &) \\
%% x_2 &=& (&\frac34&2& \frac{2k-5}4 & k-2 &) \\
%%  &\vdots \\
%% x_i &=& (&\frac{2i-1}4&i& \frac{2k-2i-1}4 & k-i &) \\
%%  &\vdots \\
%% x_{k-1} &=& (&\frac{2k-3}4&k-1& \frac14 & 1 &) \\
%% x_k &=& (&\frac{2k-1}4&k& -\frac14 & 0 &) \\
%\end{array}
%\quad\hbox{and}\quad
%\renewcommand{\arraystretch}{1.2}
%\renewcommand{\arraycolsep}{.35em}
%\begin{array}{rccclccc}
%% \\
%% y_0 &=& (&0&-1&k&2k+1&)\\
%% y_1 &=& (&\frac14&0& \frac{2k+1}4 & k+1 &)\\
%% y_2 &=& (&\frac12&1&\frac12&1&)\\
%% y_3 &=& (&\frac{2k+1}4&k+1&\frac14&0&)\\
%% y_4 &=& (&k&2k+1&0&-1&)
%% \end{array}
%% \end{equation*}
%\begin{equation*}
%\renewcommand{\arraystretch}{1.2}
%\renewcommand{\arraycolsep}{.2em}
%\begin{array}{rclccccc}
%x_0 &=& (&-\frac14&0& \frac{2k-1}4 & k &) \\
%x_1 &=& (&\frac14&1& \frac{2k-3}4 & k-1 &) \\
%x_2 &=& (&\frac34&2& \frac{2k-5}4 & k-2 &) \\
% &\vdots \\
%x_i &=& (&\frac{2i-1}4&i& \frac{2k-2i-1}4 & k-i &) \\
% &\vdots \\
%x_{k-1} &=& (&\frac{2k-3}4&k-1& \frac14 & 1 &) \\
%x_k &=& (&\frac{2k-1}4&k& -\frac14 & 0 &)
%\end{array}
%\quad\hbox{and}\quad
%\renewcommand{\arraystretch}{1.2}
%\renewcommand{\arraycolsep}{.35em}
%\begin{array}{rccclccc}
%y_0 &=& (&0&-1&k&2k+1&)\\
%y_1 &=& (&\frac14&0& \frac{2k+1}4 & k+1 &)\\
%y_2 &=& (&\frac12&1&\frac12&1&)\\
%y_3 &=& (&\frac{2k+1}4&k+1&\frac14&0&)\\
%y_4 &=& (&k&2k+1&0&-1&)
%\end{array}
%\end{equation*}

\begin{table}[ht]
$$
\renewcommand{\arraystretch}{1.2}
\renewcommand{\arraycolsep}{.2em}
\begin{array}{rclccccc}
&&&L&M&A&B\\[4pt]
x_0 &=& (&-\frac14&0& \frac{2k-1}4 & k &) \\
x_1 &=& (&\frac14&1& \frac{2k-3}4 & k-1 &) \\
x_2 &=& (&\frac34&2& \frac{2k-5}4 & k-2 &) \\
 &\vdots \\
x_i &=& (&\frac{2i-1}4&i& \frac{2k-2i-1}4 & k-i &) \\
 &\vdots \\
x_{k-1} &=& (&\frac{2k-3}4&k-1& \frac14 & 1 &) \\
x_k &=& (&\frac{2k-1}4&k& -\frac14 & 0 &) \\
%\end{array}
%\quad\hbox{and}\quad
%\renewcommand{\arraystretch}{1.2}
%\renewcommand{\arraycolsep}{.35em}
%\begin{array}{rccclccc}
\\
y_0 &=& (&0&-1&k&2k+1&)\\
y_1 &=& (&\frac14&0& \frac{2k+1}4 & k+1 &)\\
y_2 &=& (&\frac12&1&\frac12&1&)\\
y_3 &=& (&\frac{2k+1}4&k+1&\frac14&0&)\\
y_4 &=& (&k&2k+1&0&-1&)
\end{array}
$$
%\linebreak
\caption{$x_0,\dots,y_4$ in the impartial basis $L,M,A,B$.\label{tab!LMAB}}
\end{table}

\begin{pf}
Once more, either observe that these vectors satisfy all the tag
relations of the pair of long rectangles, or plug in the formulas from
\cite{BR1}, Proposition~4.1, then delete every alternate $x$ variable
(the ones tagged with a 1) and relabel to get these $x_{0\dots k}$. \end{pf}

The proof below that we can make the remaining unprojections is similar
to that of Theorem~\ref{th!22}, so we restrict ourselves to setting out the
steps and indicating how to modify them for this case.

\begin{theorem} \label{th!41}
The diptych variety on the pair of long rectangles of Figure~\ref{f!41}
exists.
\end{theorem}

\begin{pf}
First construct the 6-fold $W_0\subset
\Aff^{k+6}_{\Span{x_{0\dots k},y_2,A,B,L,M}}$ as the pullback of~$V_k$ by
the morphism \eqref{eq!subst41}. From the equations \eqref{eq!Vk} of
$V_k$, one sees that $D_0\subset W_0$, where $I_{D_0}=(x_{1\dots k},M)$,
and we can unproject this to construct $W_1$ with new ambient
variable~$y_1$.

We define $D_1\subset \Aff^{k+7}_{\Span{x_{0\dots k},y_1,y_2,A,B,L,M}}$. To
show that $D_1\subset W_1$ we check that any monomial $m$ with the same
$\T$-weight as a generator of $I_{W_1}$ that has not already been
considered is already in $I_{D_1}$. For example, if $m\simT x_iy_1$, for
any $i=1,\dots,k$, then we can suppose without loss of generality that
$m=x_0^\xi L^\la M^\mu A^\al B^\be$. By Lemma~\ref{lem!pp41}, in
impartial $L,M,A,B$ coordinates we see that
\begin{equation*}
x_iy_1 = (\tfrac i2, i, k-\tfrac i2, 2k-i+1).
\end{equation*}
His $M$-coordinate is $i\ge1$, and since $x_0=(-1/4,0,(2k-1)/4, k)$, the
only contribution to the $M$-coordinate on the right comes from $M^\mu$,
so $\mu\ge1$. In other words, $M$ divides $m$, so $m\in I_{D_1}$ as
required.

The only other equation to check has leading term $x_0y_2\simT m =
x_0^\xi y_1^\eta L^\la M^\mu A^\al B^\be$. Since both $x_0$ and $y_1$
have zero $M$ coefficient, the same argument works again. Thus
$D_1\subset W_1$, and we can unproject with new variable $y_0$ to obtain
$W_2\subset \Aff^{k+8}_{\Span{x_{0\dots k},y_{0\dots 2},A,B,L,M}}$. The
pentagrams confirm the tag equations at the bottom corners.

We continue to unproject $y_3$ and then $y_4$ to conclude. For the first
of these, define $D_2\subset \Aff^{k+8}$ by the ideal $I_{D_2} =
(x_{0\dots k-1},y_{0\dots 1},B)$ and check that $D_2\subset W_2$. We check the
critical equations (those that are not automatically in $I_{D_2}$ as a
corollary of previous checks). First suppose that $x_ky_0\simT m =
y_2^\eta L^\la M^\mu A^\al B^\be$. Since
\begin{equation*}
x_ky_0 = (\tfrac{2k-1}4,k-1,k-\tfrac14, 2k+1)
\quad
\text{and}
\quad
y_2 = (\tfrac12,1,\tfrac12,1)
\end{equation*}
consideration of the $A$-coordinate shows that $\eta < 2k$, so the
$B$-coordinate shows that $\be\ge2$; in particular, $m\in I_{D_2}$ as
required.

Now consider $y_0y_2\simT m=x_k^\xi L^\la M^\mu A^\al B^\be$. We have
\[
y_0y_2=(1/2,0,k+1/2, 2(k+1))\quad\hbox{and}\quad x_k = ((2k-1)/4,k,-1/4,0),
\]
so $\be\ge 2(k+1)$, whence $B$ divides $m$ and $m\in I_{D_2}$.

Thus we obtain $W_3\subset \Aff^{k+9}_{\Span{x_{0\dots k},y_{0\dots
3},A,B,L,M}}$ by unprojecting $D_2\subset W_2$. Finally we observe that
$D_3\subset W_3$, where $I_{D_3} = (x_{0\dots k-1},y_{0\dots 2},B)$ for
similar reasons. For example, if $y_0y_3\simT m = x_k^\xi L^\la M^\mu
A^\al B^\be$, then $y_0y_3 = (\frac{2k+1}4,k,k+1/4, 2k+1)$ and $x_k =
(\frac{2k-1}4,k,-1/4,0)$ shows that $\be\ge k+1$, so again $B$ divides
$m$ and so $m\in I_{D_3}$. Unprojecting $D_3\subset W_3$ gives the
diptych variety we seek.
\end{pf}

\subsection{Case $[1,4]$ with even $l=2k$}
Omit the even numbered $x_i$, giving Figure~\ref{f!14}. The diptych
\begin{figure}[ht]
\begin{picture}(0,110)(-50,0)
\put(45,7){\line(0,1){96}}
\put(105,7){\line(0,1){96}}
\put(45,7){\line(1,0){60}}
\put(45,103){\line(1,0){60}}
\put(45,7){\circle*5}
\put(45,22){\circle*5}
\put(45,37){\circle*5}
\put(45,73){\circle*5}
\put(45,88){\circle*5}
\put(45,103){\circle*5}
\put(15,1){$(0)$}
\put(34,49){$\vdots$}
\put(32,17){3}
\put(32,32){2}
\put(32,70){2}
\put(32,85){3}
\put(32,99){1}
\put(21,105){$A$}
\put(105,7){\circle*5}
% \put(105,31){\circle*5}
\put(105,55){\circle*5}
% \put(105,79){\circle*5}
\put(105,103){\circle*5}
\put(111,1){$(-2)$}
% \put(111,27){2}
\put(111,51){$k$}
% \put(111,75){2}
\put(111,99){2}
\put(118,105){$B$}
\end{picture}
\begin{picture}(150,110)(-220,0)
\put(45,7){\line(0,1){96}}
\put(105,7){\line(0,1){96}}
\put(45,7){\line(1,0){60}}
\put(45,103){\line(1,0){60}}
\put(45,7){\circle*5}
\put(45,22){\circle*5}
\put(45,37){\circle*5}
\put(45,73){\circle*5}
\put(45,88){\circle*5}
\put(45,103){\circle*5}
\put(32,1){1}
\put(34,49){$\vdots$}
\put(32,17){3}
\put(32,32){2}
\put(32,70){2}
\put(32,85){3}
\put(15,99){$(0)$}
\put(21,-5){$L$}
\put(105,7){\circle*5}
% \put(105,31){\circle*5}
\put(105,55){\circle*5}
% \put(105,79){\circle*5}
\put(105,103){\circle*5}
\put(111,1){2}
% \put(111,27){2}
\put(111,51){$k$}
% \put(111,75){2}
\put(111,99){$(-2)$}
\put(118,-5){$M$}
\end{picture}
\caption{Case $[4,1]$ with even $l=2k$}
\label{f!14}
\end{figure}
variety has variables $x_{0\dots k}$, $y_{0\dots 2}$, $A,B,L,M$ with the
two bottom equations
\begin{equation*}
x_1y_0= A^{2k-1}B^k+x_0L
\quad\hbox{and}\quad
x_0y_1= x_1^2A^2B+y_0^2M.
\end{equation*}
As before, adjoining $x_2,\dots,x_k$ features a long rally of flat
pentagrams, with $y_1$ against $x_{i-1},x_i$, $x_{i+1}, x_{i+2}$ and
Pfaffian equations
\begin{gather*}
y_1x_i=A^2Bx_{i+1}+L^2Mx_{i-1}, \\[4pt]
x_{i-1}x_{i+1}=x_i^2+(A^2B)^{k-i-1}(L^2M)^{i-1}AL \\[4pt]
\hbox{and}\quad
x_{i-1}x_{i+2}=x_ix_{i+1}+(A^2B)^{k-i-2}(L^2M)^{i-1}BMy_2.
\end{gather*}
These are the equations of $V_k$ after the substitution
\begin{equation*}
(a,b,c,z)\mapsto(L^2M,-y_1,A^2B,BM).
\end{equation*}

We omit the formal statement and proof of the analogue of
Theorem~\ref{th!41}: the diptych variety on the pair of long rectangles
of Figure~\ref{f!14} exists, and after the substitution the proof
unprojects $y_0$ and $y_2$ by similar arguments in impartial coordinates.

\subsection{Case $[1,4]$ with odd $l=2k+1$}
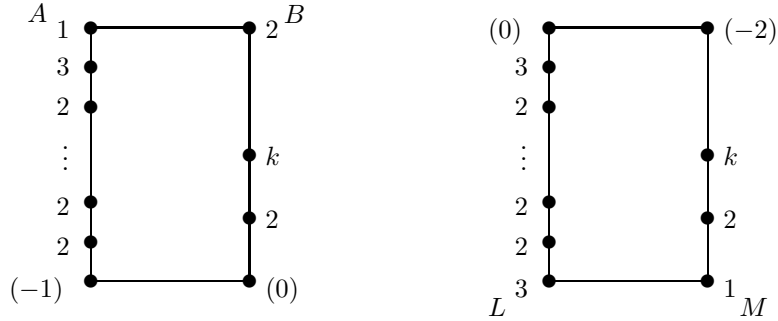
\begin{figure}[ht]
\begin{picture}(0,110)(-50,0)
\put(45,7){\line(0,1){96}}
\put(105,7){\line(0,1){96}}
\put(45,7){\line(1,0){60}}
\put(45,103){\line(1,0){60}}
\put(45,7){\circle*5}
\put(45,22){\circle*5}
\put(45,37){\circle*5}
\put(45,73){\circle*5}
\put(45,88){\circle*5}
\put(45,103){\circle*5}
\put(14,1){$(-1)$}
\put(34,49){$\vdots$}
\put(32,17){2}
\put(32,32){2}
\put(32,70){2}
\put(32,85){3}
\put(32,99){1}
\put(21,105){$A$}
\put(105,7){\circle*5}
\put(105,31){\circle*5}
\put(105,55){\circle*5}
% \put(105,79){\circle*5}
\put(105,103){\circle*5}
\put(111,1){$(0)$}
\put(111,27){2}
\put(111,51){$k$}
% \put(111,75){2}
\put(111,99){2}
\put(118,105){$B$}
\end{picture}
%\kern3em
\begin{picture}(150,110)(-220,0)
\put(45,7){\line(0,1){96}}
\put(105,7){\line(0,1){96}}
\put(45,7){\line(1,0){60}}
\put(45,103){\line(1,0){60}}
\put(45,7){\circle*5}
\put(45,22){\circle*5}
\put(45,37){\circle*5}
\put(45,73){\circle*5}
\put(45,88){\circle*5}
\put(45,103){\circle*5}
\put(22,-6){$L$}
\put(32,1){3}
\put(34,49){$\vdots$}
\put(32,17){2}
\put(32,32){2}
\put(32,70){2}
\put(32,85){3}
\put(22,99){$(0)$}
\put(105,7){\circle*5}
\put(105,31){\circle*5}
\put(105,55){\circle*5}
% \put(105,79){\circle*5}
\put(105,103){\circle*5}
\put(111,1){1}
\put(117,-6){$M$}
\put(111,27){2}
\put(111,51){$k$}
% \put(111,75){2}
\put(111,99){$(-2)$}
\end{picture}
\caption{Case $[1,4]$ with odd $l=2k+1$}
\label{f!odd}
\end{figure}

This is $[1,4]$ read from the top, but $[4,1]$ read from the bottom, so
is a mix of the two preceding cases. Omit the odd numbered $x_i$, giving
Figure~\ref{f!odd}.
The diptych variety has variables $x_{0\dots k}$, $y_{0\dots 3}$,
$A,B,L,M$ with the two bottom equations
\begin{equation*}
x_1y_0= y_1A^{2k-3}B^{k-1}+x_0^3L
\quad\hbox{and}\quad
x_0y_1= A^{2k-1}B^k+y_0M.
\end{equation*}
Adjoin $y_2$ then $x_2$ by
\begin{equation*}
\renewcommand{\arraystretch}{1.3}
\renewcommand{\arraycolsep}{.3em}
\begin{pmatrix}
y_1 & A^2B & M & y_2 \\ & y_0 & A^{2k-3}B^{k-1} & x_0^2L \\
&& x_0 & y_1 \\ &&& x_1
\end{pmatrix}
\hbox{ then }
\begin{pmatrix}
y_2 & x_1 & M & x_2 \\ & y_1 & A^2B & x_0LM \\
&& x_0 & y_2A^{2k-5}B^{k-2} \\ &&& x_1
\end{pmatrix}
\end{equation*}
After this, adjoining $x_3,\dots,x_{k-1}$ is the usual long rally of
flat pentagrams, with $y_2$ against $x_{i-1},x_i,x_{i+1},x_{i+2}$ and
\begin{equation*}
\renewcommand{\arraystretch}{1.3}
\begin{pmatrix}
y_2 & x_{i+1} & LM^2 & x_{i+2} \\
& x_{i-1} & A^2B & x_i \\
&& x_i & (A^2B)^{k-i-3}(LM^2)^{i-1}ABMy_2 \\
&&& x_{i+1}
\end{pmatrix}
\end{equation*}
and the Pfaffian equations
\begin{gather*}
y_2x_i=A^2Bx_{i+1}+LM^2x_{i-1}, \\[4pt]
x_{i-1}x_{i+1}=x_i^2+(A^2B)^{k-i-2}(LM^2)^{i-1}ABMy_2 \\[4pt]
\hbox{and}\quad
x_{i-1}x_{i+2}=x_ix_{i+1}+(A^2B)^{k-i-3}(LM^2)^{i-1}ABMy_2^2.
\end{gather*}
These are the equations of $V(k-1)$ after the substitution
\begin{equation*}
(a,b,c,z)\mapsto(LM^2,-y_2,A^2B,BM).
\end{equation*}

We again omit the formal statement and proof: the diptych variety on the
pair of long rectangles of Figure~\ref{f!odd} exists, and after the
substitution the proof unprojects $y_3$, $y_1$ and $y_0$ by arguments in
impartial coordinates.

\section{The apolar varieties $W_d$ and diptychs with $k\le3$} \label{s!lt3}

By \cite{BR1}, Classification Theorem~3.3, (3.7), when $de<3$, the cases
to treat are
\begin{equation}
\begin{matrix}
(d,e)=(1,1), & k\le2 \\
(d,e)=(1,2), & k\le3 \\
(d,e)=(1,3), & k\le5
\end{matrix}
\qquad
\begin{matrix}
(d,e)=(2,1), & k\le3 \\
(d,e)=(3,1), & k\le5
\end{matrix}
\label{eq!37}
\end{equation}
The case $k=1$ is already in \cite{BR1}, (3.9): for any values of $d,e$
we get the codimension~2 complete intersection
\begin{equation*}
\bigl(
x_1y_0 = B + Lx_0^e, \
x_0y_1 = Ax_1^d + M
\bigr)
\subset
\Aff^8_{\Span{x_0,x_1,y_0,y_1A,B,L,M}}.
\end{equation*}

In \S\ref{ss!32} we discuss the case $k=2$ for arbitrary $d,e$:
again there is an almost homogeneous variety $W_d$ that serves
as a model for the equations.

The cases with $k\ge3$ have some $x_i$ variables with tags $=1$,
which, by the tag relations \eqref{eq!xtags},
are therefore redundant generators. Eliminating them leaves a
variety in low codimension that we can specify by equations. For
$k\ge3$, the reduced models are as follows (for odd $k$, top-to-bottom
symmetry swaps $d,e$; we only list the cases with $d=1$):
\begin{equation*}
\renewcommand{\arraystretch}{1.2}
\begin{array}{c|c|c|c}
k & \text{$V_{AB}$ tags} & \text{codim as given} & \text{reduced codim}\\
\hline
3 & [1,2,1,(0)] & 4 & 2 \\
3 & [1,3,1,(0)] & 5 & 4 \\
4 & [1,3,1,3,(0)] & 5 & 4 \\
4 & [3,1,3,1,(0)] & 6 & 3 \\
5 & [1,3,1,3,1,(0)] & 6 & 2 \\
\end{array}
\end{equation*}
Eliminating the redundant generators is convenient to establish that the
varieties exist, but leaving them in has its own advantages. It allows
us to write their equations more naturally (in fact, usually as Tom
unprojections, in the language of \cite{TJ}, 2.2--2.3),
sometimes in closed Pfaffian formats. In addition, we
can put an extra deformation parameter as coefficient in front of each
variable tagged with $1$, thus exhibiting the variety as a section of a
bigger key variety.

\subsection{Case $k=2$, any $d,e$; the apolar variety $W_d$}
\label{ss!32}
For any $d,e\ge1$, the variables and tags on $V_{AB}$ are as follows:
going up the lefthand side we have $x_0,x_1,x_2$ tagged with $(0),e,d$,
against $y_{0\dots d}$ tagged with $(-e+1),2,\dots,2,1$. In $V_{AB}$ the
projection sequence first eliminates the variables $y_d,
y_{d-1},\dots,y_2$, and then the top left corner $x_2$; in $V_{LM}$ the
sequence of projections is $y_0,y_1,\dots,y_{d-2}$, then the bottom left
corner $x_0$. Following the model equations \eqref{eq!x0y0} (or \cite{BR1},~1.2),
one calculates the two equations at the bottom of the long rectangle as
\begin{equation*}
x_1y_0 = AB^d + Lx_0^d \quad\hbox{and}\quad
x_0y_1 = -x_1^{e-1} A B^{d-1} + y_0M.
\end{equation*}
One can then restore variables in the reverse order to the projection sequence
using magic pentagrams, as in \eqref{eq!pfaff}.
The $5\times5$ matrixes can be combined into a single
$(d+4)\times(d+4)$ skew matrix
\begin{equation}
\begin{pmatrix}
C & -x_0 & B & y_0 & y_1 & \dots & y_{d-1} \\
& -M & x_2 & y_1 & y_2 & \dots & y_d \\
&& x_1 & AB^{d-1} & AB^{d-2}x_2 & \dots & Ax_2^{d-1} \\
&&& Lx_0^{d-1} & LMx_0^{d-2} & \dots & LM^{d-1} \\
&&&&& \mkern-25mu \hbox{see \eqref{eq!br}} \mkern-25mu
\end{pmatrix}
\label{eq!k=2}
\end{equation}
in which we have replaced $x_1^{e-1}$ by the token $C$ in $m_{12}$; the
bottom right entries are
\begin{equation}
m_{i+5,j+5}= ALC(x_0B)^{d-j-1}(x_2M)^{i} \cdot
 \frac{(x_0x_2)^{j-i}-(BM)^{j-i}}{x_0x_2-BM}
\label{eq!br}
\end{equation}
for $0\le i<j\le d-1$.
The $4\times4$ Pfaffians of this $(d+4)\times(d+4)$ skew matrix
provide the remaining equations.

If we treat $C$ as an independent variable, then the Pfaffians of \eqref{eq!k=2} generate
the ideal of a 7-fold
\[
W_d\subset\Aff^{d+9}_{\Span{x_{0\dots 2},y_{0\dots d},A,B,L,M,C}}.
\]
It can be realised by serial unprojection following \cite{BR1},~1.2:
the equations appearing in pentagrams are
\begin{equation*}
\renewcommand{\arraystretch}{1.2}
\renewcommand{\arraycolsep}{.2em}
\begin{array}{rcl}
x_0x_2&=&-x_1C+BM \\
y_{i-1}y_{i+1}& =& y_i^2 + ALC^2(x_0B)^{d-i-1}(x_2M)^{i-1} \\
%\end{array}
%\qquad
%\renewcommand{\arraystretch}{1.2}
%\renewcommand{\arraycolsep}{.35em}
%\begin{array}{rcl}
x_0y_i &=& -x_2^{i-1}AB^{d-i}C + y_{i-1}M \\
x_1y_i &=& Ax_2^iB^{d-i} + Lx_0^{d-i}M^i \\
x_2y_i &=& y_{i+1}B - x_0^{d-i-1}CLM^i
\end{array}
\end{equation*}

The equation for $x_0x_2$ and for all $x_iy_j$ are contained among the
Pfaffians of the first 4 rows of \eqref{eq!k=2}. Beyond the 4th row,
each entry $m_{i+5,j+5}$ of \eqref{eq!br} appears in just one generating
relation,
namely
\begin{equation}
\Pf_{2,3,i+5,j+5}=Cm_{i+5,j+5}-y_iy_{j+1}+y_{i+1}y_j.
\label{eq!w2}
\end{equation}

These varieties are interesting in several ways. Replacing $x_1^{e-1}$
by the token $C$ in $m_{12}$ displays $V_{ABLM}$ as the section
$C=x_1^{e-1}$ of the 7-fold $W_d$, that is a almost homogeneous variety
under $\GL(2)\times\G_m^3$. Setting $C=0$ or $C=1$ gives invariant
6-fold sections that are also almost homogeneous. The case $d=1$ is just
the affine cone $W(1)=\aGr(2,5)$ on $\Gr(2,5)$.

\begin{exc} \rm
Write $U$ for the given representation of $\GL(2)$. Use $y_{0\dots d}$
as coefficients of a binary form $f=\sum \binom d i y_iu^{d-i}v^i \in
\Sym^dU$ and $(B,x_2)$, $(x_0,M)$ as those of two linear forms
$g=Bu+x_2v$, $h=x_0u+Mv \in U$. Then the $4\times4$ Pfaffians of
\eqref{eq!k=2} take the form
\begin{equation}
\begin{aligned}
x_1f &= Ag^d+Lh^d, \\
Mf_u-x_0f_v &= dACg^{d-1}, \\
-x_2f_u+Bf_v &= dLCh^{d-1},
\end{aligned}
\qquad
\begin{aligned}
Cx_1 & = \det\left|\begin{matrix} B & x_0 \\ x_2 & M \end{matrix}\right|
= \frac{g\wedge h}{u\wedge v}, \\
f_u \wedge f_v &= d^2ALC^2 \times \frac{ g^{d-1}\wedge h^{d-1} }{ g\wedge h },
\end{aligned}
\label{eq!Wd}
\end{equation}
where of course $f_u=\frac{\partial f}{\partial u}$ and
$f_v=\frac{\partial f}{\partial v}$. As we saw in \eqref{eq!br},
$g^{d-1}\wedge h^{d-1}$ written out as $2\times2$ minors is identically
divisible by $BM-x_0x_2$, so the final set of equations give
\eqref{eq!w2}. This form of the equations is manifestly $\GL(2)=\GL(U)$
invariant. A typical solution of \eqref{eq!Wd} is $x_0=x_2=0$ and
$x_1=A=C=L=B=M=1$, giving $g=u$, $h=v$, $f=u^d+v^d$, and one sees that
$W_d$ is the orbit closure of this typical solution under
$\GL(2)\times\G_m^3$.

At the level of the matrix \eqref{eq!k=2}, the $\GL(2)$ action replaces
rows~1 and~2 by their general linear combinations, and the $d$
rows-and-columns $5,6,\dots,d+4$ by the linear combinations
corresponding to the $(d-1)$st symmetric power. For example, adding
$\la$ times row~2 to row~1 (and the same for the columns to preserve
skew symmetry), 
\[
\la^{j-i}\times \hbox{binomial coefficient} \times \hbox{column $(5+j)$}
\]
to column $5+i$ for $j=i+1,\dots,d$ does $x_0\mapsto x_0+\la M$,
$B\mapsto B+\la x_2$ and $y_i\mapsto \sum \la^{i+j} y_j+(d-i)\la
y_{i+1}+\hbox{etc.}$, meaning $f(u,v)\mapsto f(u+\la v,v)$.
\end{exc}

\subsection{Case $k=3$; floating factors and crazy Pfaffians}
\label{ss!33}

We only need to do $e=1$; this covers $d=1$ after top-to-bottom
reflection. The case $e=1$ differs from $e\ge2$ in the order of
elimination in $V_{AB}$, as we discuss systematically in \cite{BR3}:
projecting $V_{AB}$ from the top, we eliminate $x_2$ and all the $y_i$
for $i=d-1,d-2,\dots,2$ before it becomes possible to eliminate $x_3$.
This qualitative change prevents us from treating cases with $e=1$ as a
limit of $e\ge2$.

Consider the general case $k=3$, $d\ge2$. In $V_{AB}$ we have $x_{0\dots
3}$ tagged with $(0),d,1,d$ against $y_{0,\dots,d-1}$ tagged with
$(-d+2),2,\dots,2,1$. The equations of $V_{ABLM}$ not involving $x_0$
are those of a single vertebra, and we can see them as the $4\times 4$
Pfaffians of the $(d+3)\times(d+3)$ matrix
\begin{equation}
\begin{pmatrix}
-C & x_1 & B & y_0 & y_1 & \dots & y_{d-2} \\
& LM & x_3 & y_1 & y_2 & \dots & y_{d-1} \\
&& x_2 & x_3AB^{d-2} & x_3^2AB^{d-3} & \dots & x_3^{d-1}A \\
&&& x_1^{d-2}L & x_1^{d-3}L^2M & \dots & L^{d-1}M^{d-2} \\
&&&&& \mkern-25mu \hbox{see \eqref{eq!mi5}} \mkern-25mu
\end{pmatrix}
\label{eq!m3d1}
\end{equation}
with
\begin{equation}
m_{i+5,j+5}=x_3ALC(x_1B)^{d-2-j}(x_3LM)^i\frac{(x_1x_3)^{j-i}
-(BLM)^{j-i}}{x_1x_3-BLM}.
\label{eq!mi5}
\end{equation}
For general $d$, this is the regular pullback of the apolar 7-fold
$W(d-1)$ constructed in \ref{ss!32} under the substitution
\begin{align*}
(x_{0\dots 2},y_{0\dots d-1}, & A,B,L,M,C) \\
 &\mapsto (-x_1,x_2,x_3,y_{0\dots d-1},x_3A,B,L,LM,-C).
\end{align*}
The diptych variety $V_{ABLM}$ comes from this pullback on adjoining
$x_0$ by unprojection of the divisor
\begin{align*}
D_0 &= \Aff^6_{\Span{x_1,y_0,A,B,M,C}} \\
 &= (x_2=x_3=y_{1\dots d-1}=L=0)
\subset \Aff^{d+8}_{\Span{x_{1\dots 3},y_{0\dots d-1},A,B,L,M,C}}.
\end{align*}
The Pfaffians of \eqref{eq!m3d1} clearly vanish on $D_0$, so $D_0$ is
contained in the pullback and we can unproject it to get $V_{ABLM}$.

For our application, this proves that $V_{ABLM}$ exists (for any
$d\ge2$), and we could stop there. However, this case still has a
general point to teach us: namely, how the Pfaffians of \eqref{eq!m3d1}
fit together with the unprojection equations of $x_0$.

%Starting from the other end, we see as in \cite{BR1},~1.2 that
%$V_{ABLM}$ has bottom cross
Starting from the bottom, as in \eqref{eq!x0y0}, we have
\begin{equation*}
x_1y_0 = AB^{d-1}C^2 + Lx_0 \quad\hbox{and}\quad
x_0y_1 = x_1^{d-2} AB^{d-2}C + My_0^2.
\end{equation*}
(We add a variable $C$ as annotation on $x_2$, making its tag equation
$Cx_2=x_1x_3$ in $V_{AB}$ and $V_{LM}$.) It contains the unprojection
divisor $D:(x_0=y_0=AB^{d-2}C=0)$, leading to the pentagram
$x_1,y_0,y_0,y_1,\xi$ and the $4\times4$ Pfaffians of
\begin{equation}
\begin{pmatrix}
x_1 & BC & -L & -\xi \\
& x_0 & AB^{d-2}C & -My_0 \\
&& y_0 & x_1^{d-2} \\
&&& y_1
\end{pmatrix}.
\label{eq!x1y0}
\end{equation}
The unprojection variable $\xi$ here must be $x_3$ (rather than $x_2$
with the tag $e=1$), as one sees for example from the Pfaffian
$\Pf_{12.35}=x_1^{d-1}-x_0\xi+BMCy_0$.

We link the equations together by adding a final $(d+4)$th column to
\eqref{eq!m3d1}:
\begin{equation}
\renewcommand{\arraycolsep}{0.2em}
\begin{pmatrix}
-C & x_1 & B & y_0 & y_1 & \dots & y_{d-2} & x_0 \\
& LM & x_3 & y_1 & y_2 & \dots & y_{d-1} & y_0M \\
&& x_2 & x_3AB^{d-2} & x_3^2AB^{d-3} & \dots & x_3^{d-1}A & AB^{d-1}M \\
&&& x_1^{d-2}L & x_1^{d-3}L^2M & \dots & L^{d-1}M^{d-2} & x_1^{d-1} \\
&&&&& \dots
\end{pmatrix}
\label{eq!1mc}
\end{equation}
with the same lower right entries $m_{i+5,j+5}$ as \eqref{eq!mi5}, and
the last column ending in
\begin{equation*}
m_{4+i,4+d}=-AC(Bx_1)^{d-1-i} \times \frac{(x_1x_3)^i-(BLM)^i}{x_1x_3-BLM}
\enspace\hbox{for $i=1,\dots,d-1$.}
\end{equation*}
The $4\times4$ Pfaffians of \eqref{eq!1mc} provide all but one of the
equations of $V_{ABLM}$. Comparing \eqref{eq!x1y0} with \eqref{eq!1mc},
we see that the equation
\[
x_1y_0=-AB^{d-1}C^2+x_0L
\]
is missing, although
$M$ times it is the Pfaffian $\Pf_{12.3 (d+4)}$ (in fact its multiples
by $x_1^{d-2}$, $x_2$, $x_3$, $y_1,\dots,y_{d-1}$ are also in the ideal
of Pfaffians of \eqref{eq!1mc}).

The little problem we face is how to cancel the common factor $M$ in the
entries $m_{2,3}$, $m_{2,d+4}$ and $m_{3,d+4}$ of \eqref{eq!1mc}, or in
the $3\times3$ submatrix $\left(\begin{smallmatrix} LM & y_0M \\ &
AB^{d-1}M \end{smallmatrix}\right)$ formed by rows and colums $2,3,d+4$,
without spoiling the other Pfaffians. We do this by {\em floating} $M$
from the entries with indices $2,3,d+4$ to the complementary entries
with $1,4,\dots,d+3$, adding the $4\times4$ Pfaffians of the floated
matrix, including the equation for $x_1y_0$, to those of \eqref{eq!1mc}.

The full set of equations is a mild form of {\em crazy Pfaffian},
analogous to Riemen\-schneider's quasi-determinantal \cite{R}: rather
than floating $M$ as a factor in two matrixes, we can view it as a
multiplier between entries with indices $2,3,d+4$ and those with
$1,4,\dots,d+3$; when evaluating a crazy Pfaffian, we include $M$ as a
factor whenever a product crosses between these two regions. Thus the
factors $M$ in the triangle $m_{2,3}$, $m_{2,d+4}$ and $m_{3,d+4}$ of
\eqref{eq!1mc} appear as before in most Pfaffians, but
not in $\Pf_{12.3(d+4)}$ or $\Pf_{23.i(d+4)}$ for $i=4,\dots,d+3$.

We discussed a case of floating in \cite{TJ},~9.1, especially around
(9.4), but the present instance displays the phenomenon in a
particularly clear form. This type of crazy Pfaffians or floating
factors occur frequently in our experience of working with Gorenstein
rings of codimension $\ge4$, and seem to be a basic device in
understanding how one vertebra links to the next. We expect to return to
this in future publications.

%%% In the model codimension~3 case, the Buchsbaum--Eisenbud theorem
%%% describes Pfaffian formats that provide the relations and syzygies
%%% of a Gorenstein ring exactly. Regular pullback from Gorenstein
%%% quasihomogeneous variety do the same in some cases, including most
%%% of the cases we have treated so far. But when linking one vertebra
%%% to another, we see floating factors and crazy Pfaffians almost
%%% everywhere. We expect to return to this in future publications.

%%% We treat the cases $d=2$ and $3$ here in detail, since by \eqref{eq!37}
%%% we owe those to the application. The case $d=2$ is easy: the matrixes in
%%% pentagrams each have 1 as an entry, eliminating $x_1,x_2,A$ and $L$ by
%%% equations
%%% \begin{equation*}
%%% \begin{array}{rcl}
%%% x_1 &=& x_0 x_3 - y_0 B M, \\
%%% x_2 &=& x_1 x_3 - B L M, \\
%%% \end{array}
%%% \qquad
%%% \begin{array}{rcl}
%%% A &=& x_0 y_1 - y_0^2 M, \\
%%% L &=& x_3 y_0 - y_1 B.
%%% \end{array}
%%% \end{equation*}
%%% The diptych variety is this graph over
%%% $\Aff^6_{\Span{x_0,x_3,y_0,y_1,B,M}}$.
%%%
%%% When $d=3$, \eqref{eq!1mc} becomes
%%% \begin{equation*}
%%% \begin{pmatrix}
%%% -C & x_1 & B & y_0 & y_1 & x_0 \\
%%% & LM & x_3 & y_1 & y_2 & y_0M \\
%%% && x_2 & x_3AB & x_3^2A & AB^2M \\
%%% &&& x_1L & L^2M & x_1^2 \\
%%% &&&& x_3 ALC & -x_1ABC \\
%%% &&&&& -AC(x_1x_3 + BLM)
%%% \end{pmatrix}
%%% \end{equation*}
%%% The equation giving $x_1y_0$ only appears in the ideal of these Pfaffians
%%% multiplied by $x_1$, $x_2$, $x_3$, $y_1$, $y_2$ and $M$; otherwise
%%% they give the diptych equations exactly.

\section{The cases $de=3$ and parallel unprojection}
\label{s!k45}

In \ref{s!31}, we construct all remaining cases $de=3$ with $k=4$ or 5
of \eqref{eq!37} to complete the construction of all diptych varieties with $de\le4$.
Finally, in \ref{s!key}, we observe that each of these can be realised
as a regular pullback from a single {\em key variety}, a 10-fold
$W\subset\Aff^{16}$.

\subsection{Small diptychs by pentagrams}
\label{s!31}
When $k=4$, the cases $(d,e)=(1,3)$ or $(3,1)$ are distinct. In
each case, we pass to the reduced model, which is isomorphic
to the diptych variety we seek but easier to treat because it has lower
codimension, and then adjoin the redundant generators using
pentagrams.

\begin{mycase}{$[3,1,3,1]$}
Write $x_0,x_1,x_2,x_3,x_4$ with $V_{AB}$ tags $[(0),1,3,1,3]$ opposite
$y_0,y_1,y_2$. We work up from the reduced model, that has only
$x_0,x_4$ against $y_0,y_1,y_2$; we eliminate $y_2$ from this getting
the codimension 2 complete intersection
\begin{equation*}
x_0y_1=AB+My_0 \quad\hbox{and}\quad x_4y_0=By_1^2+Lx_0,
\end{equation*}
and adjoin $y_2$ by the pentagram $x_4,x_0,y_0,y_1,y_2$ and its Pfaffian
matrix
\begin{equation*}
M_1=\begin{pmatrix}
x_4 & y_1^2 & -L & -y_2 \\
& x_0 & B & -M\\
&& y_0 & A \\
&&& y_1
\end{pmatrix}
\qquad
\begin{array}{rcl}
x_0y_2 &=& x_4A + My_1^2, \\
y_0y_2 &=& y_1^3 + AL, \\
x_4y_1 &=& y_2B + LM.
\end{array}
\end{equation*}
These five Pfaffian equations define the reduced model in codimension~3.

We recover the full set of equations by adjoining the redundant $x_2$,
then $x_1$ and $x_3$ in either order. Adjoin $x_2$ by the pentagram
$x_0,x_0,y_0$, $y_1,x_2$:
\begin{equation*}
M_2=\begin{pmatrix}
x_0 & AB & -M & -x_2 \\
& y_0 & 1 & -y_1B\\
&& y_1 & Lx_0 \\
&&& x_4
\end{pmatrix}
\qquad
\begin{array}{rcl}
x_2 &=& x_0x_4 - y_1BM\\
 &&\text{and}\\
x_2y_0 &=& y_1AB^2 + Lx_0^2,\\
x_2y_1 &=& x_4AB + LMx_0.
\end{array}
\end{equation*}
Adjoin $x_1$ by the pentagram $x_0,y_1,x_4,x_2,x_1$:
\begin{equation*}
M_3=\begin{pmatrix}
x_0 & x_2 & -BM & -x_1 \\
& y_1 & 1 & -AB\\
&& x_4 & LMx_0 \\
&&& x_2
\end{pmatrix}
\qquad
\begin{array}{rcl}
x_1 &=& x_0x_2 - AB^2M\\
 &&\text{and}\\
x_1x_4 &=& x_2^2 + x_0BLM^2,\\
x_1y_1 &=& x_2AB + LMx_0^2.
\end{array}
\end{equation*}
Finally adjoin $x_3$ by the pentagram $x_2,x_0,y_1,x_4,x_3$:
\begin{equation*}
M_4=\begin{pmatrix}
x_2 & x_4AB & -LM & -x_3 \\
& x_0 & 1 & -BM\\
&& y_1 & x_2 \\
&&& x_4
\end{pmatrix}
\qquad
\begin{array}{rcl}
x_3 &=& x_2x_4 - BLM^2\\
 &&\text{and}\\
x_0x_3 &=& x_2^2 + x_4AB^2M,\\
x_3y_1 &=& x_4^2AB + LMx_2.
\end{array}
\end{equation*}
The five Pfaffians of $M_1$ together with the three equations for
$x_1,x_2,x_3$ define $V_{ABLM}\subset\Aff^{11}_{\Span{x_{0\dots
4},y_{0\dots 1},A,B,L,M}}$.
\end{mycase}

\begin{mycase}{$[1,3,1,3]$}
Write $x_0,x_1,x_2,x_3,x_4$ with $V_{AB}$ tags $[(0),3,1,3,1]$ against
$y_0,y_1$. The reduced model is in codimension~4 on variables
$x_0,x_1,x_3,x_4,y_0,y_1$; eliminating $x_4$ then $x_3$ from this leaves
two equations
\begin{equation*}
x_0y_1 = Ax_1 + y_0^2M \quad\hbox{and}\quad x_1y_0 = A^2B + Lx_0
\end{equation*}
To recover the reduced model, we adjoin $x_3$ and then $x_4$. Adjoin
$x_3$ by the pentagram $x_1,x_0,y_0$, $y_1,x_3$:
\begin{equation*}
M_1=\begin{pmatrix}
x_1 & AB & -L & -x_3 \\
& x_0 & A & -My_0 \\
&& y_0 & x_1 \\
&&& y_1
\end{pmatrix}
\qquad
\begin{array}{rcl}
x_0x_3 &=& x_1^2 + y_0ABM, \\
x_3y_0 &=& y_1AB + x_1L, \\
x_1y_1 &=& x_3A + LMy_0.
\end{array}
\end{equation*}
The unprojection divisor of $x_4$ is $(x_0=x_1=y_0=A)$, so that the
reduced model exists. We adjoin $x_4$ by the pentagram
$x_3,x_1,y_0,y_1,x_4$:
\begin{equation*}
M_2=\begin{pmatrix}
x_3 & y_1B & -L & -x_4 \\
& x_1 & A & -LM \\
&& y_0 & x_3 \\
&&& y_1
\end{pmatrix}
\qquad
\begin{array}{rcl}
x_1x_4 &=& x_3^2 + y_1BLM, \\
x_3y_1 &=& x_4A + L^2M, \\
x_4y_0 &=& y_1^2B + x_3L.
\end{array}
\end{equation*}
These 8 equations define the reduced model in codimension~4 together
with a residual copy of $\Aff^4_{\Span{x_0,x_4,B,M}}$. Calculating with syzygies
or saturating against $y_0$ (say) recovers the {\em long equation}
\begin{equation*}
x_0 x_4 = x_1 x_3 + y_0 y_1 B M + A B L M.
\end{equation*}
In terms of the Tom and Jerry unprojections of \cite{TJ}, the
calculation to this point is a standard double Jerry; see \cite{TJ}
Section~9.2 which gives a closed form statement of the result, apart
from the long equation.

Finally, we adjoin the redundant generator $x_2$ by the pentagram
$x_1,y_0,y_1,x_3,x_2$:
\begin{equation*}
M_3=\begin{pmatrix}
x_1 & x_3A & -LM & -x_2 \\
& y_0 & 1 & -AB \\
&& y_1 & x_1L \\
&&& x_3
\end{pmatrix}
\qquad
\begin{array}{rcl}
x_1x_3 &=& x_2 + ABLM, \\
x_2y_0 &=& x_3A^2B + x_1^2L, \\
x_2y_1 &=& x_3^2A + x_1L^2M.
\end{array}
\end{equation*}
Thus the diptych in this case is the graph of $x_2=x_1x_3-ABLM$ over its
reduced model, in codimension~5 with $10\times 25$ resolution.
\end{mycase}

\begin{rem}\rm Since $x_2$ has tag 1, it makes sense to give him
annotation $C$; in the pentagram equations above, this can be done
simply by replacing the $1$ in $M_3$ by $C$. Computer algebra
experiments (after saturating these pentagram equations against $y_0LM$)
show that this gives a 7-fold $V_{ABCLM}$ in codimension 5 with
$14\times35$ resolution and serial unprojection form. (The webpage \cite{Dip}
has files to download and run in Magma \cite{Ma} to run this calculation
and other experiments.)
\end{rem}

%%% \subsection{Case $k=5$, $de=3$}
\begin{mycase}{$[1,3,1,3,1]$}
When $k=5$, we consider tags $[1,3,1,3,1]$ on $V_{AB}$; this also covers the case
$[3,1,3,1,3]$ by top-to-bottom reflection. Write
$x_0,x_1,x_2,x_3,x_4,x_5,y_0,y_1$ with $V_{AB}$ tags $[0,1,3,1,3,1]$.
The reduced model has only $x_0,x_5$ against $y_0,y_1$, with two
equations
\begin{equation*}
x_0y_1 = A + y_0M \quad\hbox{and}\quad x_5y_0 = y_1^3B + L.
\end{equation*}
The diptych variety is isomorphic to
$\Aff^6_{\Span{x_0,x_5,y_0,y_1,B,M}}$, and is the graph over it of $A,L$,
$x_4,x_2$, $x_1,x_3$ expressed as functions by
\begin{equation*}
\begin{array}{rcl}
A &=&  x_0 y_1 - y_0 M, \\
L &=& x_5 y_0 - y_1^3 B,
\end{array}
\qquad
\begin{array}{rcl}
x_1 &=& x_0 x_2 - A^2 B M, \\
x_2 &=& x_0 x_4 - y_1 A B M, \\
x_3 &=& x_2 x_4 - A B L M^2, \\
x_4 &=& x_0 x_5 - y_1^2 B M.
\end{array}
\end{equation*}
It is a fun exercise to compute all of this with magic pentagrams as in
previous cases.
\end{mycase}

\subsection{A key variety by parallel unprojection}
\label{s!key}
There is a uniform treatment of the cases $k=4$ and 5 and $de=3$ as
regular pullbacks of a key 10-fold $W$ that is given by a parallel
unprojection construction similar to that of Papadakis and Neves
\cite{PN}. We start from the codimension 2 complete intersection
$W_0\subset\Aff^{12}_{\Span{u_{1\dots4},s_{1\dots4},a_{1\dots4}}}$
given by
\begin{equation*}
\begin{aligned}
u_1u_3 &= a_2s_1s_2u_2 + a_4s_3s_4u_4, \\
u_2u_4 &= a_1s_1s_4u_1 + a_3s_2s_3u_3,
\end{aligned}
\end{equation*}
which is a normal 10-fold containing as divisors the four codimension 3
complete intersections
\begin{equation*}
(s_1,u_3,u_4), \quad (s_2,u_4,u_1), \quad (s_3,u_1,u_2), \quad (s_4,u_2,u_3).
\end{equation*}
Parallel unprojection of these four divisors gives a codimension~6
Gorenstein subvariety $W\subset\Aff^{16}_{\Span{u_{1\dots4},
v_{1\dots4}, s_{1\dots4}, a_{1\dots4}}}$ with a $20\times66$ resolution,
by standard application of the Kustin--Miller unprojection theorem. The
full set of equations is obtained as follows. Each individual
unprojection variable $v_i$ is adjoined by a pentagram, giving three
linear unprojection equations such as
\begin{equation}
\renewcommand{\arraycolsep}{0.2em}
\begin{pmatrix}
u_2 & a_1s_4u_1 & -a_3s_2s_3 & -v_1 \\
& u_3 & s_1 & -a_4s_3s_4 \\
&& u_4 & a_2s_2u_2 \\
&&& u_1
\end{pmatrix} \quad
\begin{aligned}
s_1v_1 &= u_1u_2 - a_3a_4s_2s_3^2s_4, \\
u_4v_1 &= a_1s_4u_1^2 + a_2a_3s_2^2s_3u_2, \\
u_3v_1 &= a_1a_4s_3s_4^2u_1 + a_2s_2u_2^2.
\end{aligned}
\label{eq!u1u2}
\end{equation}
In addition, there are 6 bilinear equations for $v_iv_j$, making
$2+4\times3+6=20$ equations. Four of these also come from pentagrams,
the first of which gives
\begin{equation}
v_1v_2 = a_2u_2^3 + a_1a_3a_4^2s_3^3s_4^3,
\label{eq!v1v2}
\end{equation}
whereas the remaining two are ``long equations''
\begin{equation*}
\begin{aligned}
v_1v_3 &= a_1a_4s_4^3v_4 + a_2a_3s_2^3v_2 + 3a_1a_2a_3a_4s_1s_2^2s_3s_4^2, \\
v_2v_4 &= a_1a_2s_1^3v_1 + a_3a_4s_3^3v_3 + 3a_1a_2a_3a_4s_1^2s_2s_3^2s_4
\end{aligned}
\end{equation*}
that can be computed using syzygies.

The construction has 4-fold cyclic symmetry $(1234)$, apparent in the picture
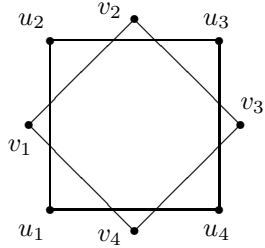
\[
\begin{picture}(80,95)(0,-8)
\put(8,8){\line(1,0){64}}
\put(8,8){\line(0,1){64}}
\put(8,8){\circle*{3}}
\put(8,72){\circle*{3}}
\put(-4,-2){$u_1$}
\put(-4,78){$u_2$}
\put(66,78){$u_3$}
\put(66,-2){$u_4$}
\put(72,72){\circle*{3}}
\put(72,8){\circle*{3}}
\put(72,72){\line(-1,0){64}}
\put(72,72){\line(0,-1){64}}
\put(40,0){\line(-1,1){40}}
\put(40,0){\line(1,1){40}}
\put(40,0){\circle*{3}}
\put(40,80){\circle*{3}}
\put(0,40){\circle*{3}}
\put(80,40){\circle*{3}}
\put(40,80){\line(-1,-1){40}}
\put(40,80){\line(1,-1){40}}
\put(26,-4){$v_4$}
\put(-8,30){$v_1$}
\put(26,82){$v_2$}
\put(80,46){$v_3$}
\end{picture}
\]
We view the $v_i$ as tagged by 1 and annotated by $s_i$ (by the first equation of
\eqref{eq!u1u2}), and the $u_i$ as tagged by 3 and annotated by $a_i$
(by \eqref{eq!v1v2}). We get Gorenstein projections on eliminating any
subset of the $v_i$, but we can only eliminate $u_i$ after projecting
out the neighbouring $v_{i-1}$ and $v_i$.

We use this variety as a model for diptych varieties.
The diptychs with $de=3$ and $k=4,5$ of \ref{s!k45} arise by pullback from
$W$ on making the following substitutions:
\begin{equation*}
\renewcommand{\arraycolsep}{0.5em}
\begin{array}{rl}
\hbox{Case } [3,1,3,1]: &
\begin{aligned}
v_1 &= x_1 \\
v_2 &= x_3 \\
v_3 &= y_2 \\
v_4 &= y_0
\end{aligned}
\qquad
\begin{aligned}
u_1 &= x_0 \\
u_2 &= x_2 \\
u_3 &= x_4 \\
u_4 &= y_1
\end{aligned}
\qquad
\begin{aligned}
a_1 &= L \\
a_2 &= 1 \\
a_3 &= A \\
a_4 &= 1
\end{aligned}
\qquad
\begin{aligned}
s_1 &= 1 \\
s_2 &= 1 \\
s_3 &= B \\
s_4 &= M
\end{aligned} \\[40pt]
\hbox{Case } [1,3,1,3]: &
\begin{aligned}
v_1 &= x_2 \\
v_2 &= x_4 \\
v_3 &= z \\
v_4 &= x_0
\end{aligned}
\qquad
\begin{aligned}
u_1 &= x_1 \\
u_2 &= x_3 \\
u_3 &= y_1 \\
u_4 &= y_0
\end{aligned}
\qquad
\begin{aligned}
a_1 &= 1 \\
a_2 &= 1 \\
a_3 &= B \\
a_4 &= M
\end{aligned}
\qquad
\begin{aligned}
s_1 &= 1 \\
s_2 &= A \\
s_3 &= 1 \\
s_4 &= L
\end{aligned} \\[40pt]
\hbox{Case } [3,1,3,1,3]: &
\begin{aligned}
v_1 &= x_1 \\
v_2 &= x_3 \\
v_3 &= x_5 \\
v_4 &= y_0
\end{aligned}
\qquad
\begin{aligned}
u_1 &= x_0 \\
u_2 &= x_2 \\
u_3 &= x_4 \\
u_4 &= y_1
\end{aligned}
\qquad
\begin{aligned}
a_1 &= L \\
a_2 &= 1 \\
a_3 &= 1 \\
a_4 &= B
\end{aligned}
\qquad
\begin{aligned}
s_1 &= 1 \\
s_2 &= 1 \\
s_3 &= A \\
s_4 &= M
\end{aligned}
\end{array}
\end{equation*}
where, in the second case, $z=y_0y_1-AL$ is a redundant generator.

\subsection*{Acknowledgments}

This work was partially supported by Korean Government WCU Grant R33-2008-000-10101-0.

\bigskip

\noindent Gavin Brown, \\
Mathematics Institute, University of Warwick, \\
Coventry CV4 7AL, England

\noindent {\it e-mail}: G.Brown@warwick.ac.uk

\bigskip

\noindent Miles Reid, \\
Mathematics Institute, University of Warwick, \\
Coventry CV4 7AL, England

\noindent {\it e-mail}: Miles.Reid@warwick.ac.uk

\end{document}